\documentclass[twoside,10pt]{siamltex}
\usepackage{amsmath,amssymb}
\usepackage{algorithm,algpseudocode}  % For pseudocode algorithms
\usepackage{listings,url,verbatim}    % Could we do without listings?
\usepackage{multicol}
\usepackage{multirow}         % For combining cells in tables
\usepackage[pdftex]{graphicx} % Laura uses this, though it restricts
\usepackage{subfigure}        % Figures in figures
\usepackage[usenames]{color}  % Laura's idea -- what's this for?

\graphicspath{{../TechReport2007/FIGURES/}}

\title{Communication-optimal parallel and sequential QR and LU factorizations}
\author{James Demmel, Laura Grigori, Mark Hoemmen, and Julien Langou}
\date{\today}

\newcommand{\RR}{{\mathbb{R}}}
\newcommand{\bmat}{\left[ \begin{array}}
\newcommand{\emat}{\end{array} \right]}

\begin{document}
% Set up program listings for C syntax
\lstset{basicstyle=\small,showstringspaces=false,language=C} 

\maketitle

\begin{abstract}
We present parallel and sequential dense QR factorization algorithms
  that are both \emph{optimal} (up to polylogarithmic factors) in the
  amount of communication they perform, and just as \emph{stable} as
  Householder QR.  
% Our first algorithm, Tall Skinny QR (TSQR), factors
% $m \times n$ matrices in a one-dimensional (1-D) block cyclic row
% layout, and is optimized for $m \gg n$.  Our second algorithm, CAQR
% (Communication-Avoiding QR), factors general rectangular matrices
% distributed in a two-dimensional block cyclic layout.  It invokes
% TSQR for each block column factorization.

% The new algorithms are superior in both theory and practice.  
We prove optimality % of sequential and parallel CAQR 
by extending known lower bounds on communication bandwidth for sequential 
and parallel matrix multiplication to provide latency lower bounds, and
show these bounds apply to the LU and QR decompositions. We not only show
that our QR algorithms attain these lower bounds (up to
polylogarithmic factors), but that existing LAPACK and ScaLAPACK
algorithms perform asymptotically more communication.
% Both TSQR and
% CAQR have asymptotically lower latency cost in the parallel case, and
% asymptotically lower latency and bandwidth costs in the sequential case.  
We also point out recent LU algorithms in the literature that
attain at least some of these lower bounds.
% In practice, we have implemented parallel TSQR on several
% machines, with speedups of up to $6.7\times$ on 16 processors of a
% Pentium III cluster, and up to $4\times$ on 32 processors of a
% BlueGene/L.  We have also implemented sequential TSQR on a laptop for
% matrices that do not fit in DRAM, so that slow memory is disk.  Our
% out-of-DRAM implementation was as little as $2\times$ slower than the
% predicted runtime as though DRAM were infinite.

% We have also modeled the performance of our parallel CAQR algorithm,
% yielding predicted speedups over ScaLAPACK's \lstinline!PDGEQRF! of up
% to $9.7\times$ on an IBM Power5, up to $22.9\times$ on a model
% Petascale machine, and up to $5.3\times$ on a model of the Grid.  
\end{abstract}

% \tableofcontents

%% mfh 07 Jul 2008: outline is in the introduction now
\section{Introduction}\label{S:introduction}

The large and increasing costs of communication motivate redesigning
algorithms to avoid it whenever possible.  
% Communication matters for both parallel and sequential algorithms.  
In the parallel case, communication refers to messages between processors, 
which may be sent over a network or via a shared memory.  
In the sequential case, communication  refers
to data movement between different levels of the memory hierarchy.
In both the parallel and sequential cases we model the
time to communicate a message of $n$ words as
$\alpha + \beta n$, where $\alpha$ is the latency and
$\beta$ is the reciprocal bandwidth.
Many authors have pointed out technology trends causing
floating point to become faster at an exponentially
higher rate than bandwidth, and bandwidth at an exponentially
higher rate than latency (see e.g., Graham et al.\ \cite{graham2005getting}).

We present parallel and sequential dense QR factorization algorithms
that are both \emph{optimal} (sometimes only up to polylogarithmic factors) 
in the amount of communication (latency and bandwidth) they require, 
and just as \emph{numerically stable} as conventional Householder QR.  
Some of the algorithms are novel, and some extend earlier work.  
The first set of algorithms, ``Tall Skinny QR'' (TSQR),
are for matrices with many more rows than columns, and the second set,
``Communication-Avoiding QR'' (CAQR), are for general rectangular
matrices.  The algorithms have significantly lower latency cost in the
parallel case, and significantly lower latency and bandwidth costs in
the sequential case, than existing algorithms in LAPACK and ScaLAPACK.
% Our algorithms are numerically stable in the same senses as in LAPACK
% and ScaLAPACK.

It will be easy to see that our parallel and sequential TSQR
implementations communicate as little as possible.
To prove optimality of CAQR, 
we extend known lower bounds on 
communication bandwidth for sequential and
parallel versions of conventional $\Theta(n^3)$ matrix multiplication 
(see Hong and Kung \cite{hong1981io}
and Irony, Toledo, and Tiskin \cite{irony2004communication}) 
to also provide latency lower bounds, and show that these 
bounds also apply to $\Theta(n^3)$ implementations of dense LU and QR 
decompositions. Showing that the bounds apply to LU is easy, 
but QR is more subtle.
We show that CAQR attains these lower bounds 
(sometimes only up to polylogarithmic factors).
% whereas existing LAPACK and ScaLAPACK algorithms perform 
% asymptotically more communication, as described further below.
% (LAPACK costs more in both latency and bandwidth, and ScaLAPACK in
% latency; it turns out that ScaLAPACK already uses optimal bandwidth.)
% Operation counts are shown in Tables
% \ref{tbl:1-par-tsqr}--\ref{tbl:6-seq-caqr-square}, and will be
% discussed below in more detail.

Implementations of TSQR and CAQR demonstrating significant speedups
over LAPACK and ScaLAPACK will be presented in other work
\cite{PRACTICE}; here we concentrate on proving optimality.

Tables~\ref{tbl:1-par-tsqr}--\ref{tbl:6-seq-caqr-square} summarize our
performance models and lower bounds for TSQR, CAQR, and 
LAPACK's sequential and ScaLAPACK's parallel QR factorizations.  
Our model of computation looks the same for the 
parallel and sequential cases, with
running time = \#flops $\times$ time\_per\_flop  + 
\#words\_moved $\times$ (1/bandwidth) + 
\#messages $\times$ latency,
% \end{eqnarray*}
where the last two terms constitute the communication.
We do not model overlap of communication and computation,
which while important in practice can at most improve the
running time by a factor of 2, whereas we are looking for
asymptotic improvements. In the tables we give the
\#flops, \#words moved and \#messages as functions of
the number of rows $m$ and columns $n$ (assuming $m \geq n$), 
the number of
processors $P$ in the parallel case, and the size of
fast memory $W$ in the sequential case.
To make these tables easier to read, we omit most lower order terms,
make boldface the terms where the new algorithms differ significantly
from Sca/LAPACK,
and make the optimal choice of matrix layout for each parallel algorithm:
This means optimally choosing the block size $b$ as well as the 
processor grid dimensions $P_r \times P_c$ in the 2-D block cyclic layout.  
(See Section~\ref{sec:CAQR_optimal}
for discussion of these parameters, and detailed performance models for
general layouts.)

Tables~\ref{tbl:1-par-tsqr}--\ref{tbl:3-par-caqr-square} present the
parallel performance models for TSQR, CAQR on general rectangular
matrices, and CAQR on square matrices, respectively.  First, Table
\ref{tbl:1-par-tsqr} shows that parallel TSQR requires only $\log P$
messages, which is both optimal and a factor $2n$ fewer messages than
ScaLAPACK's parallel QR factorization \lstinline!PDGEQRF!.  Table
\ref{tbl:2-par-caqr-general} shows that parallel CAQR needs only
$\Theta(\sqrt{nP/m})$ messages (ignoring polylogarithmic factors) on a
general $m \times n$ rectangular matrix, which is both optimal and a
factor $\Theta(\sqrt{mn/P})$ fewer messages than ScaLAPACK.  Note that
$\sqrt{mn/P}$ is the square root of each processor's local memory
size, up to a small constant factor.  Table
\ref{tbl:3-par-caqr-square} presents the same comparison for the
special case of a square $n \times n$ matrix.

Next, Tables \ref{tbl:4-seq-tsqr}--\ref{tbl:6-seq-caqr-square} 
present the sequential performance models for TSQR, CAQR on general 
rectangular matrices, and CAQR on square matrices, respectively.
Table \ref{tbl:4-seq-tsqr}
compares sequential TSQR with sequential blocked Householder QR.  This
is LAPACK's QR factorization routine \lstinline!DGEQRF! when fast
memory is cache and slow memory is DRAM, and models ScaLAPACK's
out-of-DRAM QR factorization routine \lstinline!PFDGEQRF! when fast
memory is DRAM and slow memory is disk.  Sequential TSQR transfers
fewer words between slow and fast memory: $2mn$, which is both optimal
and a factor $mn/(4W)$ fewer words than transferred by blocked
Householder QR.  Note that $mn/W$ is how many times larger the matrix
is than the fast memory size $W$.  Furthermore, TSQR requires fewer
messages: at most about $3mn/W$, which is close to optimal and $\Theta(n)$
times lower than Householder QR.
Table~\ref{tbl:5-seq-caqr-general} compares sequential CAQR and sequential
blocked Householder QR on a general rectangular matrix.  Sequential
CAQR transfers fewer words between slow and fast memory:
$\Theta(mn^2/\sqrt{W})$, which is both optimal and a factor
$\Theta(m/\sqrt{W})$ fewer words transferred than blocked Householder
QR.  Note that $m/\sqrt{W} = \sqrt{m^2/W}$ is the square root
of how many times larger a square $m \times m$ matrix is than the fast
memory size $W$.  Sequential CAQR also requires fewer messages: $12 mn^2 /
W^{3/2}$, which is optimal.  We note that our analysis of CAQR
applies for any $W$, whereas our analysis of the algorithms in
LAPACK and ScaLAPACK assume that at least 2 columns fit in fast memory,
that is $W \geq 2m$; otherwise they may communicate even more.
Finally, Table~\ref{tbl:6-seq-caqr-square} presents the same comparison for the
special case of a square $n \times n$ matrix.  

\newpage

\begin{table}[h]
\small
\centering
\begin{tabular}{l | l | l | l}
         & TSQR & \lstinline!PDGEQRF! & Lower bound \\ \hline
\# flops & $\frac{2mn^2}{P} + \frac{2n^3}{3} \log P$ 
         & $\frac{2mn^2}{P} - \frac{2n^3}{3P}$  
         & $\Theta\left( \frac{mn^2}{P} \right)$ \\ 
\# words & $\frac{n^2}{2} \log P$ 
         & $\frac{n^2}{2} \log P$ 
         & $\frac{n^2}{2} \log P$ \\
\textbf{\# messages}
         & $\mathbf{\log P}$ 
         & $\mathbf{2n \log P}$ 
         & $\mathbf{\log P}$ \\ 
\end{tabular}
\caption{Performance models of parallel TSQR and ScaLAPACK's parallel
  QR factorization \lstinline!PDGEQRF! on an $m \times n$ matrix with 
  $P$ processors, along with lower bounds on the number of flops,
  words, and messages.  We assume $m/P \geq n$.}
\label{tbl:1-par-tsqr}
\end{table}

\noindent
\rule{4.75in}{.25mm}

\begin{table}[h]
\tiny
\centering
\begin{tabular}{l | l | l | l}
         & Par.\ CAQR & \lstinline!PDGEQRF! & Lower bound \\ \hline
\# flops & $\frac{2mn^2}{P} + \frac{2n^3}{3}$
         & $\frac{2mn^2}{P} + \frac{2n^3}{3}$
         & $\Theta\left( \frac{mn^2}{P} \right)$ \\ 
\# words & $\sqrt{\frac{m n^3}{P}} \log P
            - \frac{1}{4} \sqrt{\frac{n^5}{m P}} 
              \log\left( \frac{n P}{m} \right)$    
         & $\sqrt{\frac{m n^3}{P}} \log P
            - \frac{1}{4} \sqrt{\frac{n^5}{m P}} 
              \log\left( \frac{n P}{m} \right)$    
         & $\Theta\left( \sqrt{\frac{m n^3}{P}} \right)$ \\
\textbf{\# messages} 
         & $\mathbf{ \frac{1}{4}
            \sqrt{\frac{n P}{m}}
            \log^2\left( 
                \frac{m P}{n} 
            \right) 
            \cdot \log\left( 
                P \sqrt{\frac{m P}{n}}
            \right) }$ 
         & $\mathbf{ \frac{n}{4} 
              \log\left( \frac{m P^5}{n} \right)
              \log\left( \frac{m P}{n} \right) }$
         & $\mathbf{ \Theta\left( \sqrt{\frac{nP}{m}} \right) }$ \\
\end{tabular}
\caption{Performance models of parallel CAQR and ScaLAPACK's parallel
  QR factorization \lstinline!PDGEQRF! on a $m \times n$ matrix with 
  $P$ processors, along with lower bounds on the number of flops,
  words, and messages.  The matrix is stored in a 2-D $P_r \times P_c$ 
  block cyclic layout with square $b \times b$ blocks.  We choose $b$,
  $P_r$, and $P_c$ optimally and independently for each algorithm.}
\label{tbl:2-par-caqr-general}
\end{table}

\noindent
\rule{4.75in}{.25mm}

\begin{table}[h]
\small
\centering
\begin{tabular}{l | l | l | l}
         & Par.\ CAQR & \lstinline!PDGEQRF! & Lower bound \\ \hline
\# flops & $\frac{4n^3}{3P}$
         & $\frac{4n^3}{3P}$
         & $\Theta\left( \frac{n^3}{P} \right)$ \\ 
\# words & $\frac{3n^2}{4 \sqrt{P}} \log P$
         & $\frac{3n^2}{4 \sqrt{P}} \log P$  % fixed 10 Jul 2008 
         & $\Theta\left( \frac{n^2}{\sqrt{P}} \right)$ \\
\textbf{\# messages} 
         & $\mathbf{ \frac{3}{8} \sqrt{P} \log^3 P }$
         & $\mathbf{ \frac{5n}{4} \log^2 P }$
         & $\mathbf{ \Theta\left( \sqrt{P} \right) }$ \\
\end{tabular}
\caption{Performance models of parallel CAQR and ScaLAPACK's parallel
  QR factorization \lstinline!PDGEQRF! on a square $n \times n$ matrix with 
  $P$ processors, along with lower bounds on the number of flops,
  words, and messages.  The matrix is stored in a 2-D $P_r \times P_c$ 
  block cyclic layout with square $b \times b$ blocks.  We choose $b$,
  $P_r$, and $P_c$ optimally and independently for each algorithm.}
\label{tbl:3-par-caqr-square}
\end{table}

\newpage

\begin{table}[h]
\small
\centering
\begin{tabular}{l | l | l | l}
         & Seq.\ TSQR & Householder QR & Lower bound \\ \hline
\# flops & $2mn^2$
         & $2mn^2$
         & $\Theta(mn^2)$ \\
\# words & $2mn$
         & $\frac{m^2 n^2}{2W}$
         & $2mn$ \\
\# messages & $\frac{2mn}{\widetilde{W}}$
            & $\frac{mn^2}{2W}$
            & $\frac{2mn}{W}$ \\
\end{tabular}
\caption{Performance models of sequential TSQR and blocked sequential
  Householder QR (either LAPACK's in-DRAM \lstinline!DGEQRF! or ScaLAPACK's
  out-of-DRAM \lstinline!PFDGEQRF!) on an $m \times n$ matrix with
  fast memory size $W$, along with lower bounds on the number of flops,
  words, and messages.  We assume $m \gg n$ and $W \geq 3n^2/2$.  
  $\widetilde{W} = W - n(n+1)/2$, which is at least about $\frac{2}{3}W$.}
\label{tbl:4-seq-tsqr}
\end{table}

\noindent
\rule{4.75in}{.25mm}

\begin{table}[h]
\small
\centering
\begin{tabular}{l | l | l | l}
         & Seq.\ CAQR & Householder QR & Lower bound \\ \hline
\# flops & $2mn^2 - \frac{2n^3}{3}$
         & $2mn^2 - \frac{2n^3}{3}$
         & $\Theta(mn^2)$ \\
\textbf{\# words}
         & $\mathbf{ 3 \frac{mn^2}{\sqrt{W}} }$
         & $\mathbf{ \frac{m^2 n^2}{2W} - \frac{m n^3}{6W}
            + \frac{3mn}{2} - \frac{3n^2}{4} }$ 
         & $\mathbf{ \Theta(\frac{mn^2}{\sqrt{W}} )}$ \\
%%       & $\mathbf{ \frac{3 mn^2}{32 \sqrt{2 W}} }$ \\
\textbf{\# messages}
         & $\mathbf{ 12 \frac{mn^2}{W^{3/2}} }$ 
         & $\mathbf{ \frac{mn^2}{2W} + \frac{2mn}{W} }$
         & $\mathbf{ \Theta(\frac{mn^2}{W^{3/2}} )}$ \\
%%       & $\mathbf{ \frac{3 mn^2}{32 \sqrt{2 W^3}} }$ \\
\end{tabular}
\caption{Performance models of sequential CAQR and blocked sequential
  Householder QR (either LAPACK's in-DRAM \lstinline!DGEQRF! or ScaLAPACK's
  out-of-DRAM \lstinline!PFDGEQRF!) on an $m \times n$ matrix with
  fast memory size $W$, along with lower bounds on the number of flops,
  words, and messages.}
\label{tbl:5-seq-caqr-general}
\end{table}

\noindent
\rule{4.75in}{.25mm}

\begin{table}[h]
\small
\centering
\begin{tabular}{l | l | l | l}
         & Seq.\ CAQR & Householder QR & Lower bound \\ \hline
\# flops & $\frac{4n^3}{3}$
         & $\frac{4n^3}{3}$
         & $\Theta(n^3)$ \\
\textbf{\# words}
         & $\mathbf{ 3 \frac{n^3}{\sqrt{W}} }$ 
         & $\mathbf{ \frac{n^4}{3W} + \frac{3n^2}{4} }$ 
         & $\mathbf{ \Theta(\frac{n^3}{\sqrt{W}} )}$ \\       
%%       & $\mathbf{ \frac{3n^3}{32 \sqrt{2 W}} }$ \\       
\textbf{\# messages}
         & $\mathbf{ 12 \frac{n^3}{W^{3/2}} }$ 
         & $\mathbf{ \frac{n^3}{2W} }$
         & $\mathbf{ \Theta(\frac{n^3}{W^{3/2}} )}$ \\       
%%       & $\mathbf{ \frac{3n^3}{32 \sqrt{2 W^3}} }$ \\       
\end{tabular}
\caption{Performance models of sequential CAQR and blocked sequential
  Householder QR (either LAPACK's in-DRAM \lstinline!DGEQRF! or ScaLAPACK's
  out-of-DRAM \lstinline!PFDGEQRF!) on a square $n \times n$ matrix with
  fast memory size $W$, along with lower bounds on the number of flops,
  words, and messages.}
\label{tbl:6-seq-caqr-square}
\end{table}

%% mfh 29 Jul 2008
%\newpage

Finally, we note that although our new algorithms perform
slightly more floating point operations than LAPACK and 
ScaLAPACK, they have the same highest order terms in
their floating point operation counts.
(For TSQR, which is intended for the case $m \gg n$,
only the term containing $m$ is highest order.)
In fact we prove a matching lower bound on the amount of arithmetic,
assuming we avoid ``Strassen-like'' algorithms in a way
made formal later.

\begin{comment}
We have concentrated on the cases of a homogeneous parallel computer
and a sequential computer with a two-level memory hierarchy. But real
computers are obviously more complicated, combining many levels of
parallelism and memory hierarchy, perhaps heterogeneously. 
We partially address this more difficult problem in two ways.
First, we show that our parallel and sequential TSQR designs correspond 
to the two simplest cases of reduction trees (binary and flat, respectively),
and that different choices of reduction trees will let us optimize
TSQR for more general architectures.
Second, we describe how to apply TSQR and CAQR recursively
to accommodate hierarchical architectures, analogously to
multiple levels of blocking for matrix multiplication.
\end{comment}

Now we briefly describe related work and our contributions.
The tree-based QR idea itself is not novel (see for example,
\cite{buttari2007class,buttari2007parallel,cunha2002new,golub1988parallel,gunter2005parallel,kurzak2008qr,pothen1989distributed,quintana-orti2008scheduling,rabani2001outcore}),
but we have a number of optimizations and generalizations:
\begin{itemize}
% \item Some authors (e.g., \cite{golub1988parallel}) formulate the
%   inner steps of the factorization to require $\Theta(n \log P)$
%   communication steps, a factor of $\Theta(n)$ larger than the optimal
%   number of steps (which our factorization achieves).
\item Our algorithm can perform almost all its floating-point
  operations using any fast sequential QR factorization routine.  
  For example, we can use blocked Householder transformation
  exploiting BLAS3, or invoke Elmroth
  and Gustavson's recursive QR (see
  \cite{elmroth1998new,elmroth2000applying}).

\item We use TSQR as a building block for CAQR, for the parallel
  resp.\ sequential factorization of arbitrary rectangular matrices in
  a two-dimensional block cyclic layout.

\item Most significantly, we prove optimality for both our parallel 
  and sequential algorithms, with a 1-D layout for TSQR and 2-D block 
  layout for CAQR, i.e., that they minimize bandwidth and latency costs.  
  This assumes $\Theta(n^3)$ (non-Strassen-like) algorithms, and is usually 
  shown in a Big-Oh sense, sometimes modulo polylogarithmic terms.

\item We describe special cases in which existing sequential algorithms 
by Elmroth and Gustavson \cite{elmroth2000applying} and also LAPACK's DGEQRF
attain minimum bandwidth. In particular, with the correct choice of
block size, Elmroth's and Gustavson's RGEQRF algorithm attains minimum
bandwidth and flop count, though not minimum latency.

\item We observe that there are alternative LU algorithms in 
  the literature that attain at least some of these communication
  lower bounds: \cite{grigori2008calu} describes a parallel LU algorithm
  attaining both bandwidth and latency lower bounds, and
  \cite{toledo1997locality} describes a sequential LU algorithm that
  at least attains the bandwidth lower bound.

\item 
  We outline how to extend both algorithms and optimality results
  to certain kinds of hierarchical architectures, either with multiple
  levels of memory hierarchy, or multiple levels of parallelism
  (e.g., where each node in a parallel machine consists of other parallel
   machines, such as multicore).
   In the case of TSQR we do this by adapting it to work on general 
   reduction trees.  
\end{itemize}

It is possible to do a stable QR factorization (or indeed most any
dense linear algebra operation) at the same asymptotic speed as
matrix multiplication (e.g., in $\Theta(n^{\log_2 7})$ operations using
Strassen) \cite{FastLinearAlgebraIsStable} and so with less
communication as well, but we do not discuss these algorithms in
this paper.

We note that the $Q$ factor will be represented as a tree of smaller $Q$
factors, which differs from the traditional layout.  Many previous
authors did not explain in detail how to apply a stored TSQR $Q$
factor, quite possibly because this is not required for solving 
a single least squares problem:
Adjoining the right-hand side(s) to the matrix $A$,
and taking the QR factorization of the result, requires only the $R$
factor.  Previous authors discuss this optimization.  However, many of
our applications require storing and working with the implicit
representation of the $Q$ factor.  
Our performance models show that applying this tree-structured $Q$
has about the same cost as the traditionally represented $Q$.
% Furthermore, applying this implicit
% representation has nearly the same performance model as constructing an
% explicit $Q$ factor with the same dimensions as $A$.

% \subsection{Outline}\label{SS:intro:outline}

The rest of this report is organized as follows.  
Section~\ref{sec:TSQR_optimal} presents TSQR,
describing its parallel and sequential optimizations,
performance models, comparisons to LAPACK and ScaLAPACK,
and how it can be adapted to other architectures.
Section~\ref{sec:CAQR_optimal} presents CAQR analogously.
(This paper is based on the technical report
\cite{TSQR_technical_report}, to which we leave many of the detailed
derivations of the performance models.)
Section~\ref{sec:LowerBounds_TSQR} presents our lower bounds
for TSQR, and
Section~\ref{sec:LowerBounds_CAQR} for CAQR (as well as LU).
Section~\ref{S:related-work} describes related work.
Section~\ref{sec:Conclusions_optimal} summarizes
and describes open problems and future work.

\section{Tall-Skinny QR - TSQR}
\label{sec:TSQR_optimal}

In this section, we present the TSQR algorithm for
computing the QR factorization of an $m \times n$ matrix $A$, 
stored in a 1-D block row layout.
We assume $m \geq n$, and typically $m \gg n$.
(See \cite{scalapackusersguide} for a description of 1D and 2D layouts.)

Subsection~\ref{sec:TSQR_optimal_tree} describes 
parallel TSQR on a binary tree, 
sequential TSQR on a ``flat'' tree, and then TSQR
as a reduction on an arbitrary tree.
Subsection~\ref{sec:TSQR_optimal_perfmodel} describes 
performance models, 
and Subsection~\ref{sec:TSQR_optimal_comparison} 
compares TSQR to  alternative algorithms,
both stable and unstable;
we will see that TSQR does asymptotically less
communication than the stable alternatives, and is
about as fast as the fastest unstable alternative
(CholeskyQR).
% Subsection~\ref{sec:TSQR_optimal_relatedwork} describes 
% related work.

\subsection{TSQR as a reduction operation}
\label{sec:TSQR_optimal_tree}

We will describe a family of algorithms that takes
an $m$-by-$n$ matrix $A = [A_0 ; A_1 ; \cdots ; A_{p-1}]$
and produces the $R$ factor of its QR decomposition.
Here we use Matlab notation, so that the $A_i$ are
stacked atop one another, and we assume $A_i$ is 
$m_i$-by-$n$. In later sections we will assume
$m_i \geq n$, but that is not necessary here.
% with $m_i \geq n$. 
% {\em (Do we need $m_i \geq n$ for correctness?
% What if we choose all $Q$ matrices square?
% Doesn't it work for general $m_i$?)}
% Our illustrations below will use $p=4$ blocks.

The basic operation in our examples is to take
two or more matrices stacked atop one another,
like $\hat{A} = [A_0; A_1]$, and replace them by 
the $R$ factor of $\hat{A}$. As long as 
more than one matrix remains in the stack, the
reduction continues until one $R$ factor is left,
which we claim is the $R$ factor of the original $A$.
The pattern of which pairs (or larger groups) of matrices 
are combined in one step forms what we will call a reduction tree.

We write this out explicitly for TSQR performed
on a binary tree starting with $p=4$ blocks.
We start by replacing each $A_i$ by its own
individual $R$ factor:

\begin{equation}
\label{eqn:TSQR_binarytree_1}
A = 
\begin{pmatrix}
A_0 \\
A_1 \\
A_2 \\
A_3 \\
\end{pmatrix}
=
\begin{pmatrix}
  Q_{0} R_{0} \\
  Q_{1} R_{1} \\
  Q_{2} R_{2} \\
  Q_{3} R_{3} \\
\end{pmatrix}.
\end{equation}

Proceeding with the first set of reductions, we write
\begin{equation}
\label{eqn:TSQR_binarytree_2}
\begin{pmatrix}
  R_{0} \\
  R_{1} \\ \hline
  R_{2} \\ 
  R_{3} \\
\end{pmatrix}
= 
\begin{pmatrix}
  \begin{pmatrix} 
    R_{0} \\
    R_{1} \\
  \end{pmatrix} \\ \hline
  \begin{pmatrix}
    R_{2} \\
    R_{3} \\
  \end{pmatrix}
\end{pmatrix}
=
\begin{pmatrix}
  Q_{01} R_{01} \\ \hline
  Q_{23} R_{23} \\ 
\end{pmatrix}
\end{equation}
Thus $[R_0;R_1]$ is replaced by $R_{01}$
and $[R_2;R_3]$ is replaced by $R_{23}$.
Here and later, the subscripts on a matrix like $R_{ij}$ refer to
the original $A_i$ and $A_j$ on which they depend.

The next and last reduction is
\begin{equation}
\label{eqn:TSQR_binarytree_3}
\begin{pmatrix}
  R_{01} \\
  R_{23} \\
\end{pmatrix}
=
Q_{0123} R_{0123}.
\end{equation}
We claim that $R_{0123}$ is the $R$ factor 
of the original $A=[A_0;A_1;A_2;A_3]$.
To see this, we combine 
equations~(\ref{eqn:TSQR_binarytree_1}),
(\ref{eqn:TSQR_binarytree_2}) and
(\ref{eqn:TSQR_binarytree_3}) to write

\begin{equation}
\label{eqn:TSQR_binarytree_4}
A  
= 
\begin{pmatrix}
A_0 \\
A_1 \\
A_2 \\
A_3 \\
\end{pmatrix}
= 
\left(
  \begin{array}{c | c | c | c}
    Q_{0} & & & \\ \hline
    & Q_{1} & & \\ \hline
    & & Q_{2} & \\ \hline
    & & & Q_{3} \\ 
  \end{array}
\right)
\cdot
\left(
  \begin{array}{c | c}
    Q_{01} &        \\ \hline
          &  Q_{23} \\ 
  \end{array}
\right)
\cdot
Q_{0123} \cdot R_{0123}
\end{equation}

\begin{comment}
\begin{eqnarray*}
A  
& = &
\begin{pmatrix}
A_0 \\
A_1 \\
A_2 \\
A_3 \\
\end{pmatrix}
= 
\left(
  \begin{array}{c | c | c | c}
    Q_{0} & & & \\ \hline
    & Q_{1} & & \\ \hline
    & & Q_{2} & \\ \hline
    & & & Q_{3} \\ 
  \end{array}
\right)
\cdot
\begin{pmatrix}
R_0 \\
R_1 \\
R_2 \\
R_3 \\
\end{pmatrix}
\\
& = & 
\left(
  \begin{array}{c | c | c | c}
    Q_{0} & & & \\ \hline
    & Q_{1} & & \\ \hline
    & & Q_{2} & \\ \hline
    & & & Q_{3} \\ 
  \end{array}
\right)
\cdot
\left(
  \begin{array}{c | c}
    Q_{01} &        \\ \hline
          &  Q_{23} \\ 
  \end{array}
\right)
\cdot
\begin{pmatrix}
  R_{01} \\ \hline
  R_{23} \\
\end{pmatrix} \\
& = &
\left(
  \begin{array}{c | c | c | c}
    Q_{0} & & & \\ \hline
    & Q_{1} & & \\ \hline
    & & Q_{2} & \\ \hline
    & & & Q_{3} \\ 
  \end{array}
\right)
\cdot
\left(
  \begin{array}{c | c}
    Q_{01} &        \\ \hline
          &  Q_{23} \\ 
  \end{array}
\right)
\cdot
Q_{0123} \cdot R_{0123}
\; \; .
\end{eqnarray*}
\end{comment}

For this product to make sense, we must 
choose the dimensions of the $Q$ factors consistently: 
They can all be square, 
or when all $m_i \geq n$, 
they can all have $n$ columns (in which case each $R$ factor 
will be $n$-by-$n$). (The usual representation of $Q$ factors by 
Householder vectors encodes both possibilities.)
In either case, we have expressed $A$ as a product of
(block diagonal) orthogonal matrices
(which must therefore also be orthogonal), and the triangular
matrix $R_{0123}$. By uniqueness of the
QR decomposition (modulo signs of diagonal
entries of $R_{0123}$), this is the QR decomposition
of $A$. We note that we will not multiply the various
$Q$ factors together, but leave them represented by the
``tree of $Q$ factors'' implied by
equation~(\ref{eqn:TSQR_binarytree_4}).

We abbreviate this algorithm with the following simple
notation, which makes the binary tree apparent:

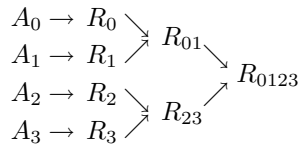
\begin{center}
\setlength{\unitlength}{.5cm}
\begin{picture}(7,4)

\put(0.5,0.5){$A_3$}
\put(0.5,1.5){$A_2$}
\put(0.5,2.5){$A_1$}
\put(0.5,3.5){$A_0$}

\put(1.5,0.5){$\rightarrow$}
\put(1.5,1.5){$\rightarrow$}
\put(1.5,2.5){$\rightarrow$}
\put(1.5,3.5){$\rightarrow$}

\put(2.5,0.5){$R_3$}
\put(2.5,1.5){$R_2$}
\put(2.5,2.5){$R_1$}
\put(2.5,3.5){$R_0$}

\put(3.5,0.65){$\nearrow$}
\put(3.5,1.35){$\searrow$}
\put(3.5,2.65){$\nearrow$}
\put(3.5,3.35){$\searrow$}

\put(4.5,1.0){$R_{23}$}
\put(4.5,3.0){$R_{01}$}

\put(5.6,1.5){$\nearrow$}
\put(5.6,2.5){$\searrow$}

\put(6.5,2.0){$R_{0123}$}

\end{picture}
\end{center}

The notation has the following meaning: if one or more
arrows point to the same matrix, that matrix
is the $R$ factor of the matrix obtained by stacking
all the matrices at the other ends of the arrows atop
one another.
This notation not only makes the parallelism in the algorithm
apparent (all QR decompositions at the same depth in the
tree can potentially be done in parallel), but implies that 
{\em any} tree leads to a valid QR decomposition. For example, conventional
QR decomposition may be expressed as the trivial tree

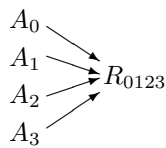
\begin{center}
\setlength{\unitlength}{.5cm}
\begin{picture}(4,4)

\put(0.5,0.5){$A_3$}
\put(0.5,1.5){$A_2$}
\put(0.5,2.5){$A_1$}
\put(0.5,3.5){$A_0$}

\put(1.5,0.85){\vector(3,2){1.4}}
\put(1.5,1.70){\vector(3,1){1.4}}
\put(1.5,2.75){\vector(3,-1){1.4}}
\put(1.5,3.55){\vector(3,-2){1.4}}

\put(3.0,2.0){$R_{0123}$}

\end{picture}
\end{center}

The tree we will use for sequential TSQR with limited
fast memory $W$ is the following so-called ``flat tree'':

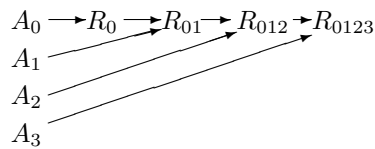
\begin{center}
\setlength{\unitlength}{.5cm}
\begin{picture}(7,4)

\put(0.5,0.5){$A_3$}
\put(0.5,1.5){$A_2$}
\put(0.5,2.5){$A_1$}
\put(0.5,3.5){$A_0$}

\put(1.5,1.0){\vector(3,1){7}}
\put(1.5,1.75){\vector(3,1){5}}
\put(1.5,2.75){\vector(4,1){3}}
\put(1.5,3.75){\vector(1,0){1}}

\put(2.5,3.5){$R_0$}
\put(3.5,3.75){\vector(1,0){1}}
\put(4.5,3.5){$R_{01}$}
\put(5.55,3.75){\vector(1,0){.8}}
\put(6.5,3.5){$R_{012}$}
\put(8.0,3.75){\vector(1,0){.5}}
\put(8.5,3.5){$R_{0123}$}

\end{picture}
\end{center}

The idea of sequential TSQR is that if fast memory can only hold a little more than a fraction $m/p$ 
of the rows of $A$ (a little more than $m/4$ for the above tree), then the algorithm proceeds
by reading in the first $m/p$ rows of $A$, doing its QR decomposition,
keeping $R_0$ in fast memory but writing the representation of $Q_0$
back to slow memory,
and then repeatedly reading in the next $m/p$ rows, doing the QR decomposition 
of them stacked below the $R$ factor already in memory,
and writing out the representation of the new $Q$ factor.
This way the entire matrix is read into fast memory once, and
the representation of all the $Q$ factors is written out to fast memory once, 
which is clearly the minimal amount of data movement possible.

For an example of yet another TSQR reduction tree more suitable for 
a hybrid parallel / out-of-core factorization, see 
\cite[Section 4.3]{TSQR_technical_report}.

It is evident that all these variants of TSQR are numerically stable,
since they just involve repeated applications of orthogonal
transformations.  Note also that the local QR factorizations in both
the parallel and sequential TSQR algorithms can avoid storing and
performing arithmetic with zeros in the triangular factors.  This
optimization still allows the use of high-performance QR algorithms
(such as the BLAS 3 $Y T Y^T$ representation of Schreiber and Van Loan
\cite{schreiber1989storage} and the recursive QR factorization of
Elmroth and Gustavson \cite{elmroth2000applying}) for the local
computations.  For details, see Demmel et al.\ \cite[Section
7]{TSQR_technical_report}.

We close this subsection by observing that the general theory of
reduction operations applied to associative operators
(e.g., optimizing the shape of the reduction tree
\cite{nishtala2008performance}, or how to compute prefix sums 
of $a_1 \star a_2 \star \cdots \star a_p$
where $\star$ could be scalar addition, matrix multiplication, etc.) 
applies to QR decomposition as well, because the mapping from
$[A_0;A_1]$ to its $R$ factor is associative (modulo roundoff and
the choice of the signs of the diagonal entries).

\subsection{Performance models for TSQR}
\label{sec:TSQR_optimal_perfmodel}

We present performance models for parallel and sequential TSQR.
We outline their derivations, which are straightforward based
on the previous descriptions, and leave details to
\cite[Section 8]{TSQR_technical_report}.
In the next section we will compare the models for TSQR
with alternative algorithms.

The runtimes will be functions of $m$ and $n$.
In the parallel case, the runtime will also depend on
the number of processors $P$, where we assume each
processor stores $m/P$ rows of the input matrix $A$.
(It is easiest to think of the rows as contiguous,
but if they are not, we simply get the QR decomposition
of a row-permutation of $A$, which is still just the QR
decomposition). In the sequential case the runtime will
depend on $W$, the size of fast memory.
We assume fast memory is large enough to 
contain at least $n$ rows of $A$, and an $R$ factor,
i.e. $W \stackrel{>}{\approx} \frac{3}{2}n^2$.
In both parallel and sequential cases, we let
$\gamma = $ time per flop, $\beta = $ reciprocal
bandwidth (time per word) and $\alpha = $ latency
(time per message). We assume no overlap of
communication and computation (as said before, this
could speed up the algorithm at most 2$\times$).
All logarithms are in base 2.

A parallel TSQR factorization on a binary reduction tree performs the
following computations along the critical path: one local QR
factorization of a fully dense $m/P \times n$ matrix, and $\log P$
factorizations, each of a $2n \times n$ matrix consisting of two $n
\times n$ upper triangular matrices.  The factorization requires
$\frac{2mn^2}{P} + \frac{2n^3}{3} \log P$
flops (ignoring lower order terms here and elsewhere) and 
$\log P$ messages, and transfers a total of 
$\frac{1}{2} n^2 \log P$ 
words between processors.  
Thus, the total run time is
\begin{equation}
\label{eqn:TSQR_par_runtime}
\text{Time}_{\text{Par.\ TSQR}}(m,n,P) = 
\left(
  \frac{2mn^2}{P} + \frac{2n^3}{3} \log P
\right) \gamma + 
\left(
  \frac{1}{2} n^2 \log P
\right) \beta + 
\left( \log P \right) \alpha 
\; \; .
\end{equation}

Now we consider sequential TSQR.
To first order, TSQR performs the same number of floating point operations
as standard Householder QR, namely $2mn^3 - 2n^3 / 3$.
As described before, sequential TSQR moves $2mn$ words by dividing 
$A$ into submatrices that are as large as possible, i.e., $m'$ rows
each such that $m' \cdot n + \frac{n(n+2)}{2} \leq W$,
or $m' \approx (W - \frac{n(n+1)}{2})/n  = \widetilde{W}/n$,
where $\widetilde{W} = W - \frac{n(n+1)}{2}$. Assuming
$A$ is stored so that groups of $m'$ rows are in contiguous
memory locations, the number of messages sequential TSQR needs
to send is $\frac{2mn}{m' n} = \frac{2mn}{\widetilde{W}}$.
Thus the runtime for sequential TSQR is
\begin{equation}
\label{eqn:TSQR_seq_runtime}
\text{Time}_{\text{Seq.\ TSQR}}(m,n,W) = 
\left(
  2mn^2 - \frac{2n^3}{3}
\right) \gamma + 
\left( 2mn \right) \beta + 
\left(
  \frac{2mn}{\widetilde{W}}
\right) \alpha 
\; \; .
\end{equation}
We note that $\widetilde{W} \stackrel{>}{\approx} 2 W / 3$, so that
the number of messages
$2mn / \widetilde{W} \stackrel{<}{\approx} 3mn / W$.

\subsection{Comparison of TSQR to alternative algorithms}
\label{sec:TSQR_optimal_comparison} 

We compare parallel and sequential QR to alternative algorithms, both
stable and unstable: Classical Gram-Schmidt (CGS), Modified Gram-Schmidt (MGS),
Cholesky QR, and Householder QR, as implemented in LAPACK and ScaLAPACK;
only the latter are numerically stable in all cases.
In summary, TSQR not only has the lowest complexity (comparing highest order terms), 
but has asymptotically lower communication complexity than the only
numerically stable alternatives.
We outline our approach and leave details of counting to
\cite[Section 9]{TSQR_technical_report}.

MGS and CGS can be either right-looking or left-looking. For CGS either
alternative has the same communication complexity, but for MGS the right-looking
variant has much less latency, so we present its performance model.

Cholesky QR forms $A^TA$, computes its upper triangular Cholesky
factor $R$, and forms $Q = A R^{-1}$. It can obviously be unstable, but
is frequently used when $A$ is expected to be well-conditioned
(see section~\ref{S:related-work}).

\begin{comment}
In contrast to parallel TSQR, parallel Householder QR needs $\log P$ messages
for each of the $n$ columns in order to compute the $n$
Householder vectors. A straightforward analysis of the
right-looking ScaLAPACK algorithm PDGEQRF
with the same 1D layout as TSQR (the left-looking version is similar)
yields a runtime of
\begin{equation}
\label{eqn:1D_ScaLAPACK_par_runtime}
Time_{1D ScaLAPACK}(m,n,P) = 
(\frac{2mn^2}{P} - \frac{2n^3}{3})  \gamma + 
(\frac{1}{2} n^2 \log P) \beta + 
(2n \log P) \alpha 
\; \; .
\end{equation}
We see that while TSQR does more flops than ScaLAPACK, it is a lower order term
(since we assume $m \gg n$). 
{\em (Need to fix Tables 1 and/or 10 in tech report.)}
TSQR has the same bandwidth requirement as
ScaLAPACK. Most importantly, TSQR sends $2n$ times fewer messages.
\end{comment}

We need to say a little more about sequential Householder QR.
LAPACK's right-looking DGEQRF repeatedly sweeps over the entire
matrix, potentially leading to proportionally as much memory traffic 
as there are floating point operations, a factor $\Theta (n)$
more than sequential TSQR; a left-looking version of 
DGEQRF would be similar.
% {\em (Begs the question of Elmroth-Gustavson again.)}
To make a fairer comparison, we model the performance
of a left-looking QR algorithm that was optimized to
minimize memory movement in an out-of-DRAM environment,
i.e., where fast memory is DRAM and slow memory is disk.
This routine, PFDGEQRF \cite{dazevedo1997design} was
designed to combine ScaLAPACK's parallelism with minimal
disk accesses.
As originally formulated, it uses ScaLAPACK's
parallel QR factorization PDGEQRF to perform the current
panel factorization in DRAM,
but we assume here that it is running sequentially since we are
only interested in modeling the traffic between slow and fast memory.
PFDGEQRF is a left-looking method, as usual with out-of-DRAM algorithms
(left-looking schemes do fewer writes than right-looking schemes,
since writes are often more expensive.)
PFDGEQRF keeps two panels in memory: a left
panel of fixed width $b$, and the current panel being factored, whose
width $c$ can expand to fill the available memory.  
Details of the algorithm and analysis may be found in
\cite{dazevedo1997design} and
\cite[Appendix F]{TSQR_technical_report}, where we choose
$b$ and $c$ to minimize disk traffic;
we summarize the performance model in Table~\ref{tbl:TSQR:perfcomp:seq}.

\begin{comment}
\begin{eqnarray}
\label{eqn:1D_PFDGEQRF_seq_runtime}
Time_{1D seq. PFDGEQRF}(m,n,W) 
& = &
(2mn^2 - \frac{2n^3}{3}) \gamma 
\nonumber \\
& &
+ 
(\frac{mn^2}{2W}(m - \frac{n}{3}) +\frac{3n}{2}(m - \frac{n}{2})) \beta  
\nonumber \\ 
& &
+ (\frac{mn}{2W}(n+4)) \alpha 
\; \; .
\end{eqnarray}
In other words, it performs the same number of floating point operations
as sequential TSQR (or DGEQRF), but moves about $\frac{mn}{4W}$ times as
many words. Note that $\frac{mn}{W}$ is how many times larger the matrix
is than fast memory. It also sends about $\frac{n}{4}$ times as many
messages as sequential TSQR.
\end{comment}

\begin{table}
  \centering
  \begin{tabular}{l|c|c|c}
    Parallel algorithm & \# flops & \# messages & \# words \\ \hline
    TSQR & $\frac{2mn^2}{P} + \frac{2n^3}{3} \log(P)$
         & $\log(P)$
         & $\frac{n^2}{2} \log(P)$ \\
    \lstinline!PDGEQRF!
        & $\frac{2mn^2}{P} - \frac{2n^3}{3P}$ 
% FIXME (mfh 24 Jun 2008): what about blocking?
        & $2n \log(P)$
        & $\frac{n^2}{2} \log(P)$ \\
% FIXME above here
    MGS    & $\frac{2mn^2}{P}$
          & $2n \log(P)$
          & $\frac{n^2}{2}\log(P)$ \\
    CGS    & $\frac{2mn^2}{P}$
          & $2n \log(P)$ 
          & $\frac{n^2}{2}\log(P)$ \\ 
%     CGS2 & $\frac{4mn^2}{P}$ 
%          & $4n \log(P)$ 
%          & $n^2 \log(P)$ \\
    CholeskyQR & $\frac{2mn^2}{P} + \frac{n^3}{3}$ 
               & $\log(P)$ 
               & $\frac{n^2}{2}\log(P)$ \\ 
  \end{tabular}
  \caption{Performance models of various parallel QR 
    algorithms for "tall-skinny" matrices, i.e. with $m \gg n$.
    We show only the best-performing versions of MGS (right-looking) 
    and CGS (left-looking).}
%   ``CGS2'' means CGS with one reorthogonalization pass.
%   Lower-order terms omitted.  All parallel terms are counted along
%   the critical path.  
%   We omit CGS2 because it is no slower than applying 
%   CGS twice, but the number of orthogonalization steps may vary 
%   based on the numerical properties of the input, so it is hard to 
%   predict performance \emph{a priori}.}
  \label{tbl:TSQR:perfcomp:par}
\end{table}

\begin{table}
\small
  \centering
  \begin{tabular}{l|c|c|c}
    Sequential algorithm & \# flops & \# messages & \# words  \\\hline
    TSQR & $2mn^2 - \frac{2n^3}{3}$ 
         & $\frac{2mn}{\widetilde{W}}$
         & $2mn - \frac{n(n+1)}{2} 
           + \frac{mn^2}{\widetilde{W}}$ \\
    \lstinline!PFDGEQRF!
         & $2mn^2 - \frac{2 n^3}{3}$
         & $\frac{2mn}{W} + \frac{mn^2}{2W}$
         & $\frac{m^2 n^2}{2W} - \frac{mn^3}{6W}
            + \frac{3mn}{2} - \frac{3n^2}{4}$ \\
    MGS & $2mn^2$ 
        & $\frac{2mn^2}{\widetilde{W}}$
        & $\frac{3mn}{2} + \frac{m^2 n^2}{2 \widetilde{W}}$ \\
    CholeskyQR & $2mn^2 + \frac{n^3}{3}$
               & $\frac{6mn}{W}$ 
               & $3mn$ \\
  \end{tabular}
  \caption{Performance models of various sequential QR 
    algorithms for "tall-skinny" matrices, i.e. with $m \gg n$.
    \lstinline!PFDGEQRF! is our model of ScaLAPACK's
    out-of-DRAM QR factorization; $W$ is the fast memory size, and
    $\widetilde{W} = W - n(n+1)/2$.  Lower-order terms omitted.}
%   We omit CGS2 because it is no slower than applying CGS twice, but the
%   number of orthogonalization steps may vary based on the numerical
%   properties of the input, so it is hard to predict performance
%   \emph{a priori}.}
  \label{tbl:TSQR:perfcomp:seq}
\end{table}

Examining Table~\ref{tbl:TSQR:perfcomp:par}, we see that
all parallel algorithms have the same highest order term in their
flop counts, $\frac{2mn^2}{P}$, and also use the same bandwidth, $\frac{n^2}{2} \log P$,
but that parallel TSQR sends
$2n$ times fewer messages than the only stable alternative
(PDGEQRF), and is about as fast as the fastest unstable method
(Cholesky QR). In other words, only parallel TSQR is simultaneously
fastest and stable.

Examining Table~\ref{tbl:TSQR:perfcomp:seq}, we see a similar
story, with sequential TSQR sending about $\frac{mn}{4W}$ times fewer
words and $\frac{n}{4}$ times fewer messages than the only
stable alternative, PFDGEQRF. Note that $\frac{mn}{W}$ is how
many times larger the entire matrix is than fast memory.
Since we assume $W \geq n^2$, the number of words TSQR sends is
less than the number of words CholeskyQR sends.

%%%%%%%%%%%%%%%%%%%%%%%%%%%%%%%%%%%%%%%%%%%%%%%%%%%%%%%%%%%%%%%%%%%%%%
%%%%%%%%%%%%%%%%%%%%%%%%%%%%%%%%%%%%%%%%%%%%%%%%%%%%%%%%%%%%%%%%%%%%%%

\section{Communication-Avoiding QR - CAQR}
\label{sec:CAQR_optimal}

We present the CAQR algorithm for computing the
QR factorization of an $m$-by-$n$ matrix $A$, with $m \geq n$.
In the parallel case $A$ is stored on a two-dimensional grid of 
processors $P = P_r \times P_c$ in a 2-D block-cyclic layout,
with blocks of dimension $b \times b$.  
We assume that all the blocks have the same size; 
we can always pad the input matrix with zero rows and columns to
ensure this is possible.  
In the sequential case we also assume $A$ is stored in a
$P_r \times P_c$ 2-D blocked layout, with individual 
$\frac{m}{P_r}$-by-$\frac{n}{P_c}$ 
blocks stored contiguously in memory.
For a detailed description of the 2-D 
block cyclic layout, see \cite{scalapackusersguide}.

Stated most simply, parallel (resp. sequential) CAQR simply
implements the right-looking QR factorization using parallel 
(resp. sequential) TSQR as the panel factorization.
The rest is bookkeeping. 

Section~\ref{sec:CAQR_parallel} discusses parallel CAQR in
more detail, 
% including choosing the parameters $P_r$, $P_c$ 
% and $b$ to optimize the runtime, 
and comparing performance to ScaLAPACK. 
We also show, given $m$, $n$ and $P$, to
choose $P_r$, $P_c$ and $b$ to minimize running times
of both algorithms; our proof of CAQR's optimality depends 
on these choices.  Section~\ref{sec:CAQR_sequential} does the
same for sequential CAQR and an out-of-DRAM algorithm from
ScaLAPACK, whose floating point operations are counted
sequentially.
Subsection~\ref{sec:seq_qr_other} discusses other sequential
QR algorithms, including showing that recursive QR routines of Elmroth
and Gustavson \cite{elmroth2000applying} also minimize
bandwidth, though possibly not latency.

\subsection{Parallel CAQR}
\label{sec:CAQR_parallel}

We describe a few details
most relevant to the complexity but refer the reader to
\cite[Section 13]{TSQR_technical_report} for details.
At the $j$-th step of the algorithm, parallel TSQR
is used to factor the panel of dimension $m-(j-1)b$-by-$b$,
whose top left corner is at matrix diagonal entry $(j-1)b+1$.
We assume for simplicity that the $m_j = m-(j-1)b$ rows are
distributed across all $P_r$ processors in the processor column.
When we do parallel TSQR on the panel, all the at most
$\frac{m}{P_r}$ local rows of the panel stored on a processor are factored together 
in the first step of TSQR. After the panel factorization,
we multiply the transpose of the $Q$ factor times the trailing
submatrix as follows. First, the Householder vectors representing the $Q$ factor 
of the $\frac{m}{P_r}$ local rows of the panel are broadcast to all the processes in the
same processor row, and applied to their submatrices in an
embarrassingly parallel fashion.
Second, the Householder vectors $Y$ of the smaller $Q$ factors in TSQR's
binary reduction tree are independently broadcast along their processor rows, 
and the updates to the $b$ rows in each pair of processors are performed
in parallel, with the triangular $T$ factor of the block Householder
transformation $I - YTY^T$ being computed by one of the two processors,
and with the two processors exchanging only $b$ rows of data.

Table~\ref{tbl:CAQR:par:model} summarizes the operation counts,
including divisions counted separately, as well as a similar
model for ScaLAPACK's PDGEQRF for comparison.
We make the following observations.
Parallel CAQR does slightly more
flops than ScaLAPACK (but only in lower order terms), and sends
nearly the same of words (actually very slightly fewer). 
But CAQR reduces the $3n \log P_r$ term in ScaLAPACK's message
count by a factor of $b$, and so can reduce the overall
message count by as much as a factor of $b$ (depending $P_r$ and $P_c$).
Thus by increasing the block size $b$, we can lower the number of messages
by a large factor. But we can't raise $b$ arbitrarily without
increasing the flop count; next we show how to
choose the parameters $b$, $P_r$ and $P_c$ to minimize the runtime.

\begin{table}[h]
\small
\centering
\begin{tabular}{l | l}
            & Parallel CAQR \\ \hline
\# messages & $\frac{3n}{b} \log P_r + \frac{2n}{b} \log P_c$ \\ \hline
\# words    & $\left( 
                   \frac{n^2}{P_c} 
                   + \frac{bn}{2} 
               \right) \log P_r
               + \left( 
                   \frac{mn - n^2/2}{P_r} + 2n 
               \right) \log P_c$ \\ \hline
\# flops    & $\frac{2n^2(3m-n)}{3P} 
               + \frac{bn^2}{2P_c} 
               + \frac{3bn(2m - n)}{2P_r} 
               + \left( \frac{4 b^2 n}{3} 
                   + \frac{n^2 (3b+5)}{2 P_c} 
               \right) \log P_r
               - b^2 n$ \\ \hline
\# divisions & $\frac{mn - n^2/2}{P_r} 
                + \frac{bn}{2} \left( \log P_r - 1 \right)$
             \\ \hline \hline
             & ScaLAPACK's \texttt{PDGEQRF} \\ \hline
\# messages & $3n \log P_r + \frac{2n}{b} \log P_c$ \\ \hline
\# words    & $\left( 
                   \frac{n^2}{P_c} 
                   + bn 
               \right) \log P_r
               + \left( 
                   \frac{mn - n^2/2}{P_r} 
                   + \frac{bn}{2} 
               \right) \log P_c$ \\ \hline
\# flops    & $\frac{2n^2(3m-n)}{3P} 
               + \frac{bn^2}{2P_c} 
               + \frac{3bn(2m - n)}{2P_r}
               - \frac{b^2 n}{3 P_r}$
            \\ \hline
\# divisions & $\frac{mn - n^2/2}{P_r}$ \\
\end{tabular}
\caption{Performance models of parallel CAQR and ScaLAPACK's
  \lstinline!PDGEQRF! when factoring an $m \times n$ matrix, $m \geq n$,
  distributed in a 2-D block cyclic layout on a $P_r \times P_c$ 
  grid of processors with square $b \times b$ blocks.  All terms are
  counted along the critical path.  In this table exclusively, ``flops'' 
  only includes floating-point additions and multiplications, not 
  floating-point divisions, which are shown separately.  
  Some lower-order terms are omitted.}
% Note that the number of flops,
% divisions, and words transferred all roughly match between the two
% algorithms, but the number of messages is about $b$ times lower 
% for CAQR.}
\label{tbl:CAQR:par:model}
\end{table}

When choosing $b$, $P_r$, and $P_c$ to minimize the runtime,
they must satisfy the following conditions:
\begin{equation}\label{eq:CAQR:par:opt:ansatz:constraints}
1 \leq P_r, P_c \leq P  
\; \; , \; 
P_r \cdot P_c = P      
\; \; , \; 
1 \leq b \leq \frac{m}{P_r} 
\; \; {\rm and} \; 
1 \leq b \leq \frac{n}{P_c} 
\end{equation}
For simplicity we will assume that $P_r$ evenly divides $m$ and that $P_c$
evenly divides $n$.  
Example values of $b$, $P_r$, and $P_c$ which satisfy the
constraints in Equation \eqref{eq:CAQR:par:opt:ansatz:constraints} are
\[
P_r = \sqrt{\frac{m P}{n}} 
\; \; , \; 
P_c = \sqrt{\frac{n P}{m}} 
\; \; {\rm and} \; 
b   = \sqrt{\frac{m n}{P}} 
\]
These values are chosen simultaneously to minimize the approximate
number of words sent, $n^2/P_c + mn/P_r$, and the approximate number
of messages, $5n/b$, where for simplicity we temporarily ignore
logarithmic factors and lower-order terms in Table
\ref{tbl:CAQR:par:model}.  This suggests using the following ansatz:
\begin{equation}\label{eq:CAQR:par:opt:ansatz}
P_r = K \cdot \sqrt{\frac{m P}{n}}
\; \; , \;
P_c = \frac{1}{K} \cdot \sqrt{\frac{n P}{m}}
\; \; \text{and} \;
b   = B \cdot \sqrt{\frac{m n}{P}}, 
\end{equation}
for general values of $K$ and $B \leq \min\{ K, 1/K \}$, since we can
thereby explore all possible values of $b$, $P_r$ and $P_c$ satisfying
\eqref{eq:CAQR:par:opt:ansatz:constraints}.

Using the substitutions in Equation \eqref{eq:CAQR:par:opt:ansatz},
the flop count (neglecting lower-order terms, including the division
counts) becomes
\begin{multline}\label{eq:CAQR:par:opt:ansatz:flops}
\frac{mn^2}{P} \left( 
    2 
    - B^2 
    + \frac{3B}{K} 
    + \frac{B K}{2} 
\right)
- 
\frac{n^3}{P} \left( 
    \frac{2}{3} 
    + \frac{3B}{2K} 
\right)
+ \\
\frac{mn^2 \log\left( K \cdot \sqrt{\frac{mP}{n}} \right)}{P} \left( 
    \frac{4B^2}{3} 
    + \frac{3 B K}{2}
\right).
\end{multline}
We wish to choose $B$ and $K$ so as to minimize the flop count.  We
know at least that we need to eliminate the dominant $mn^2
\log(\dots)$ term, so that parallel CAQR has the same asymptotic flop
count as ScaLAPACK's \lstinline!PDGEQRF!.  This is because we know
that CAQR performs at least as many floating-point operations
(asymptotically) as \lstinline!PDGEQRF!, so matching the highest-order
terms will help minimize CAQR's flop count.

To make the high-order terms of \eqref{eq:CAQR:par:opt:ansatz:flops}
match the $2mn^2/P - 2n^3/(3P)$ flop count of ScaLAPACK's parallel QR
routine, while minimizing communication as well, we can pick $K=1$ and 
\[
B = o\left(
    \log^{-1}\left(
        \sqrt{\frac{ m P }{ n }}
    \right)
\right);
\]
for simplicity we will use 
\begin{equation}\label{eq:CAQR:par:opt:flops:B}
B = \log^{-2} \left(
    \sqrt{\frac{ m P }{ n }}
\right)
\end{equation}
although $B$ could be multiplied by 
some positive constant.

The above choices of $B$ and $K$ make the flop count as follows, with
some lower-order terms omitted:
\begin{equation}\label{eq:CAQR:par:opt:flops}
\frac{2mn^2}{P}
- \frac{2n^3}{3P}
+ \frac{3 m n^2}{P \log\left( \frac{m P}{n} \right)}
\end{equation}
Thus, we can choose the block size $b$ so as to match the higher-order
terms of the flop count of ScaLAPACK's parallel QR factorization
\lstinline!PDGEQRF!.

Using the substitutions in Equations \eqref{eq:CAQR:par:opt:ansatz}
and \eqref{eq:CAQR:par:opt:flops:B} with $K = 1$, the number of
messages becomes
\begin{equation}
%\label{eq:CAQR:par:opt:lat:C}
\label{eq:CAQR:par:opt:lat}
% \frac{1}{C}
\sqrt{\frac{n P}{m}}
\cdot \log^2\left( \sqrt{\frac{m P}{n}} \right)
\cdot \log\left( P \sqrt{\frac{m P}{n}} \right).
\end{equation}

\begin{comment}
The best we can do with the latency is to make $C$ as large as
possible, which makes the block size $b$ as large as possible.  The
value $C$ must be a constant, however; specifically, the flop counts
require
\[
\begin{aligned}
C &= \Omega\left( \log^{-2} \left( K \sqrt{\frac{m P}{n}} \right)
\right)\,\text{and} \\
C &= O(1).
\end{aligned}
\]
We leave $C$ as a tuning parameter in the number of messages Equation
\eqref{eq:CAQR:par:opt:lat}.
\end{comment}

Using the substitutions in Equation \eqref{eq:CAQR:par:opt:ansatz}
and \eqref{eq:CAQR:par:opt:flops:B}, the number of words transferred
between processors on the critical path, neglecting lower-order terms,
becomes
% fixed 10 Jul 2008
\begin{multline}\label{eq:CAQR:par:opt:bw}
\sqrt{\frac{m n^3}{P}} \log P 
- \frac{1}{4} \sqrt{\frac{n^5}{m P}} \log \left( \frac{n P}{m} \right)
+ \frac{1}{4} \sqrt{\frac{m n}{P}} \log^3\left( \frac{m P}{n} \right)
\approx \\
\sqrt{\frac{m n^3}{P}} \log P 
- \frac{1}{4} \sqrt{\frac{n^5}{m P}} \log \left( \frac{n P}{m} \right).
\end{multline}

\begin{comment}
In the second step above, we eliminated the $C$ term, as it is a
lower-order term (since $m \geq n$).  Thus, $C$ only has a significant
effect on the number of messages and not the number of words
transferred.
\end{comment}

The results of these computations are shown in
Table \ref{tbl:CAQR:par:model:opt}, which also shows 
the results for ScaLAPACK, whose analogous
analysis appears in
\cite[Section 15]{TSQR_technical_report},
and the communication lower bounds, which are discussed
in Section~\ref{sec:LowerBounds_CAQR}.

\begin{comment}
the number of messages and
number of words used by parallel CAQR and ScaLAPACK when $P_r$, $P_c$,
and $b$ are independently chosen so as to minimize the runtime models,
as well as the optimal choices of these parameters.  In summary, if we
choose $b$, $P_r$, and $P_c$ independently and optimally for both
algorithms, the two algorithms match in the number of flops and words
transferred, but CAQR sends a factor of $\Theta(\sqrt{mn/P})$ messages
fewer than ScaLAPACK QR.  This factor is the local memory requirement
on each processor, up to a small constant.
\end{comment}

\begin{table}[h!]
\centering
\begin{tabular}{l | l}
             & Parallel CAQR w/ optimal $b$, $P_r$, $P_c$ \\ \hline
\# flops     & $\frac{2mn^2}{P} - \frac{2n^3}{3P}$ \\
\# messages  & $\frac{1}{4}
                \sqrt{\frac{n P}{m}}
                \log^2\left( 
                    \frac{m P}{n} 
                \right) 
                \cdot \log\left( 
                    P \sqrt{\frac{m P}{n}}
                \right)$ \\
\# words & $\sqrt{\frac{m n^3}{P}} \log P
            - \frac{1}{4} \sqrt{\frac{n^5}{m P}} 
              \log\left( \frac{n P}{m} \right)$ \\
%% Note (mfh 27 Jun 2008) -- it's 4C and not just C because the
%% original formula had log^{-2}(\sqrt{mP/n})$, and the square root
%% becomes 1/2 which is then squared (1/4) and inverted (4).  But C is
%% just a constant so it eats the 4.
Optimal $b$    & $ \sqrt{\frac{m n}{P}} 
                  \log^{-2} \left( \frac{m P}{n} \right)$ \\
Optimal $P_r$  & $\sqrt{\frac{m P}{n}}$ \\
Optimal $P_c$  & $\sqrt{\frac{n P}{m}}$ 
               \\ \hline \hline
%% PDGEQRF counts fixed (mfh 10 Jul 2008)
             & \texttt{PDGEQRF} w/ optimal $b$, $P_r$, $P_c$ \\ \hline
\# flops     & $\frac{2mn^2}{P} - \frac{2n^3}{3P}$ \\
\# messages & $\frac{n}{4} 
                 \log\left( \frac{m P^5}{n} \right)
                 \log\left( \frac{m P}{n} \right)
               + \frac{3n}{2} \log\left( \frac{m P}{n} \right)$ \\
\# words & $\sqrt{\frac{m n^3}{P}} \log P
            - \frac{1}{4} \sqrt{\frac{n^5}{m P}} 
              \log\left( \frac{n P}{m} \right)$ \\
Optimal $b$    & $ \sqrt{\frac{mn}{P}} 
                  \log^{-1}\left( \frac{m P}{n} \right)$ \\
Optimal $P_r$  & $\sqrt{\frac{m P}{n}}$ \\
Optimal $P_c$  & $\sqrt{\frac{n P}{m}}$
               \\ \hline \hline
             & Theoretical lower bound \\ \hline
%% change 2^13 to 2^11 in following two lines
\# messages  & $\sqrt{\frac{n P}{2^{11} m}}$ \\ 
\# words     & $\sqrt{\frac{mn^3}{2^{11} P}}$ \\
\end{tabular}
\caption{Highest-order terms in the performance models of parallel
	CAQR, ScaLAPACK's \lstinline!PDGEQRF!, and theoretical lower bounds
  for each, when factoring an $m \times n$ matrix, distributed in a
  2-D block cyclic layout on a $P_r \times P_c$ grid of processors
  with square $b \times b$ blocks.  All terms are counted along the
  critical path.  The theoretical lower bounds assume that $n \geq
  2^{11} m / P$, i.e., that the matrix is not too tall and skinny.
% The parameter $C$ in both algorithms is a $\Theta(1)$ tuning parameter.  
  In summary, if we choose $b$, $P_r$, and $P_c$ independently and
optimally for both algorithms, the two algorithms match in the number of flops and words
transferred, but CAQR sends a factor of $\Theta(\sqrt{mn/P})$ messages
fewer than ScaLAPACK QR.
This factor is the local memory requirement on each processor, up to a
small constant.}
\label{tbl:CAQR:par:model:opt}
\end{table}

\subsection{Sequential CAQR}
\label{sec:CAQR_sequential}

As stated above, sequential CAQR is just right-looking QR factorization
with TSQR used for the panel factorization. (In fact left-looking
QR with TSQR has the same costs \cite[Appendix C]{TSQR_technical_report},
but we stick with the right-looking algorithm for simplicity.)
We also assume the $m$-by-$n$ matrix $A$ is stored in a
$P_r \times P_c$ 2-D blocked layout, with individual 
$\frac{m}{P_r}$-by-$\frac{n}{P_c}$ 
blocks stored contiguously in memory, with $m \geq n$ and
$\frac{m}{P_r} \geq \frac{n}{P_c}$. 

For TSQR to work as analyzed we need to choose
$P_r$ and $P_c$ large enough for one such 
$\frac{m}{P_r}$-by-$\frac{n}{P_c}$ 
block to fit in fast memory, plus a bit more.
For CAQR we will need to choose $P_r$ and $P_c$ a bit
larger, so that a bit more than 3 such blocks fit in fast memory;
this is in order to perform an update on two such blocks
in the trailing matrix given Householder vectors from TSQR
occupying $\frac{mn}{P_r P_c} + \frac{n^2}{2P_c^2}$ words,
or at most $\frac{4mn}{P}$ altogether. In other words,
we need $\frac{4mn}{P}\leq W$ or $P \geq \frac{4mn}{W}$.

Leaving details to \cite[Appendix C]{TSQR_technical_report},
we summarize the complexity analysis by

\begin{eqnarray} 
\label{eq:CAQR:seq:modeltime:P}
T_{\text{seq.\ CAQR}} (m,n,P_c,P_r) 
& \leq & \left( \frac{3}{2} P(P_c-1) \right) \alpha +  \nonumber \\
&      & \left( \frac{3}{2} mn \left( P_c + \frac{4}{3} \right) 
                   - \frac{1}{2} n^2 P_c 
%          + O\left(n^2 + nP\right) 
         \right) \beta \\
&      & + \left( 2n^2m - \frac{2}{3}n^3 \right) \gamma  \nonumber
\end{eqnarray}
where we have ignored lower order terms, 
and used $P_r$ as an upper
bound on the number of blocks in each panel
since this only increases the run time slightly, 
and is simpler to evaluate than for the true number of
blocks $P_r - \lfloor (J-1) \frac{n P_r}{m P_c} \rfloor$.

Now we choose $P$, $P_r$ and $P_c$ to minimize the runtime.
From the above formula for
$T_{\text{seq.\ CAQR}} (m,n,P_c,P_r)$, we see that the runtime is
an increasing function of $P_r$ and $P_c$, so that we would like
to choose them as small as possible, within the limits imposed by
the fast memory size $P \geq \frac{4mn}{W}$. So we choose
$P = \frac{4mn}{W}$ (assuming here and elsewhere that the
denominator evenly divides the numerator). But we still need to
choose $P_r$ and $P_c$ subject to $P_r \cdot P_c = P$.

Examining $T_{\text{seq.\ CAQR}} (m,n,P_c,P_r)$ again, we see
that if $P$ is fixed, the runtime is also an increasing function
of $P_c$, which we therefore want to minimize. But we are assuming
$\frac{m}{P_r} \geq \frac{n}{P_c}$, or $P_c \geq \frac{nP_r}{m}$.
The optimal choice is therefore  $P_c = \frac{nP_r}{m}$
or $P_c = \sqrt{\frac{nP}{m}}$, which also means
$\frac{m}{P_r} = \frac{n}{P_c}$, i.e., the blocks in the algorithm
are square. This choice of $P_r = \frac{2m}{\sqrt{W}}$ and 
$P_c = \frac{2n}{\sqrt{W}}$ therefore minimizes the runtime,
yielding
\begin{eqnarray} 
\label{eq:CAQR:seq:modeltime:P:opt}
T_{\text{Seq.\ CAQR}} (m,n,W) & \leq &
      \left( 12 \frac{mn^2}{W^{3/2}} \right) \alpha +
      \left( 3 \frac{mn^2}{\sqrt{W}} +
%     + O\left( \frac{mn^2}{W} \right) 
     \right) \beta + 
     \nonumber \\
&   &      \left( 2mn^2 - \frac{2}{3}n^3 \right) \gamma.
\end{eqnarray}

We note that the bandwidth term is proportional to $\frac{mn^2}{\sqrt{W}}$,
and the latency term is $W$ times smaller,
both of which match (to within constant factors), the lower bounds
on bandwidth and latency to be described in
Section~\ref{sec:LowerBounds_CAQR}.

The results of this analysis are shown in 
Table~\ref{tbl:CAQR:seq:model:opt}, which also
shows the results for an out-of-DRAM algorithm
PFDGEQRF from ScaLAPACK, whose internal block
sizes $b$ and $c$ have been chosen to minimize disk traffic,
and where we count the floating point operations sequentially
(see \cite[Appendix F]{TSQR_technical_report});
it can also be thought of as a hypothetical model 
for an optimized left-looking version of LAPACK's DGEQRF.

%% FIXED with optimal parameter choices (mfh 23 Jun 2008)
\begin{table}[h]
\centering
\begin{tabular}{l | l}
             & Sequential CAQR w/ optimal $P_c$, $P_c$ \\ \hline
\# flops     & $2mn^2 - \frac{2}{3}n^3$  \\
\# messages  & $12 \frac{mn^2}{W^{3/2}}$ \\
\# words     & $3 \frac{mn^2}{\sqrt{W}}$ \\
Opt.\ $P$    & $4mn/W$ \\
Opt.\ $P_r$  & $2m / \sqrt{W}$ \\
Opt.\ $P_c$  & $2n / \sqrt{W}$
               \\ \hline \hline
             & ScaLAPACK's \texttt{PFDGEQRF} w/ optimal $b$, $c$\\ \hline
\# flops     & $2mn^2 - \frac{2}{3}n^3$  \\
\# messages  & $\frac{mn^2}{2W} + \frac{2mn}{W}$ \\
\# words     & $\frac{m^2 n^2}{2W} - \frac{m n^3}{6W}
                + \frac{3mn}{2} - \frac{3n^2}{4}$ \\
Opt.\ $b$    & $1$ \\
Opt.\ $c$    & $\approx \frac{W}{m}$
               \\ \hline \hline
             & Theoretical lower bound \\ \hline
% \# messages  & $\frac{\frac{mn^2}{4} - \frac{n^2}{8} \left( 
%                    \frac{n}{2} + 1 \right)}{\sqrt{8W^3}}$ \\
\# messages  & $\frac{3n^2(m - \frac{4}{3})}{16(8W^3)^{1/2}} - 1$ \\
% \# words     & $\frac{\frac{mn^2}{4} - \frac{n^2}{8} \left( 
%                    \frac{n}{2} + 1 \right)}{\sqrt{8W}}$ \\
\# words     & $\frac{3n^2(m - \frac{4}{3})}{16(8W)^{1/2}} - W$ \\

\end{tabular}
\caption{Highest-order terms in the performance models of sequential
  CAQR, ScaLAPACK's out-of-DRAM QR factorization \texttt{PFDGEQRF} 
  running on one processor, and theoretical lower bounds for each, when
  factoring an $m \times n$ matrix with a fast memory capacity of $W$ words.}
% In the case of sequential CAQR, the matrix is arranged in a
% 2-D block layout on a $P_r \times P_c$ grid of $P$ blocks.  (See 
% Appendix \ref{SS:CAQR-seq-detailed:opt}).  The optimal choices of
% these parameters result in square blocks (i.e., $m/P_r = n/P_c$).  
% In the case of \texttt{PFDGEQRF}, the $b$ parameter is the left 
% panel width and the $c$ parameter is the current panel width.  (See
% Appendix \ref{S:PFDGEQRF}.)}
\label{tbl:CAQR:seq:model:opt}
\end{table}

\subsection{Other Bandwidth Minimizing Sequential QR Algorithms}
\label{sec:seq_qr_other}

In this section we describe special cases in which previous
sequential QR algorithms also minimize bandwidth, although
they do not minimize latency.
In particular, we discuss 
two variants of Elmroth's and Gustavson's recursive
QR (RGEQR3 and RGEQRF \cite{elmroth2000applying}),
as well as LAPACK's DGEQRF.

The fully recursive routine RGEQR3 is analogous to Toledo's
fully recursive LU routine \cite{toledo1997locality}: Both
routines factor the left half of the matrix (recursively), 
use the resulting factorization of the left half to update 
the right half, and then factor the right half (recursively again).
The base case consists of a single column. The output of
RGEQR3 applied to an $m$-by-$n$ matrix returns the $Q$ 
factor in the form $I-YTY^T$, where $Y$ is the $m$-by-$n$
lower triangular matrix of Householder vectors,
and $T$ is an $n$-by-$n$ upper triangular matrix.
A simple recurrence for the number of memory references
of either RGEQR3 or Toledo's algorithm is
\begin{eqnarray}
\label{eqn:RGEQR3}
B(m,n) & = & \left\{ \begin{array}{ll} 
             B(m,\frac{n}{2}) + B(m-\frac{n}{2},\frac{n}{2}) + 
             O(\frac{mn^2}{\sqrt{W}})
             & {\rm if} \; mn > W \; {\rm and} \; n>1 \\
             mn & {\rm if} \; mn \leq W  \\
             m  & {\rm if} \; m > W \; {\rm and} \; n=1 
             \end{array} \right. \nonumber \\
    & \leq & \left\{ \begin{array}{ll} 
             2B(m,\frac{n}{2}) + 
             O(\frac{mn^2}{\sqrt{W}})
             & {\rm if} \; mn > W \; {\rm and} \; n>1 \\
             mn & {\rm if} \; mn \leq W  \\
             m  & {\rm if} \; m > W \; {\rm and} \; n=1 
             \end{array} \right. \nonumber \\
       & = & O(\frac{mn^2}{\sqrt{W}}) + mn
\end{eqnarray}
% B(m,n) & = & B(m,\frac{n}{2}) + B(m-\frac{n}{2},\frac{n}{2}) +
%              O(\max(mn, \frac{mn^2}{\sqrt{W}})) \nonumber \\
%        & \leq & 2B(m,\frac{n}{2}) + 
%              O(\max(mn, \frac{mn^2}{\sqrt{W}})) \nonumber \\
%        & = & O(\frac{mn^2}{\sqrt{W}} + mn \log W) \\
%        & = & O(\frac{mn^2}{\sqrt{W}}) 
%            \; \; {\rm if} \; \sqrt{W} \log W = O(n) \nonumber
So RGEQR3 attains our bandwidth lower bound.
(The $mn$ term must be included to account for the case
when $n<\sqrt{W}$, since each of the $mn$ matrix entries
must be accessed at least once.)
However, RGEQR3 does 
a factor greater than one
% about $\frac{7}{6}$? 3/2? 6?
times as many floating point operations
as sequential Householder QR.

Now we consider RGEQRF and DGEQRF, which are both
right-looking algorithms and differ only in how
they perform the panel factorization (by RGEQR3
and DGEQR2, resp.). Let $b$ be the width of
the panel in either algorithm. It is easy to
see that a reasonable estimate of the number of 
memory references just for the updates by all the panels
is the number of panels $\frac{n}{b}$ times the minimum
number of memory references for the average
size update $\Theta(\max(mn,\frac{mnb}{\sqrt{W}}))$,
or $\Theta(\max(\frac{mn^2}{b},\frac{mn^2}{\sqrt{W}}))$.
Thus we need to pick $b$ at least about as large
as $\sqrt{W}$ to attain the desired lower bound
$O(\frac{mn^2}{\sqrt{W}})$.

Concentrating now on RGEQRF, we get from 
inequality~(\ref{eqn:RGEQR3})
that the $\frac{n}{b}$ panel factorizations using RGEQR3
cost at most an additional \linebreak
$O(\frac{n}{b} \cdot [\frac{mb^2}{\sqrt{W}} + mb] )
= O( \frac{mnb}{\sqrt{W}} + mn)$
memory references, or $O(mn)$ if we pick $b=\sqrt{W}$.
Thus the total number of memory references for RGEQRF
with $b= \sqrt{W}$ is $O(\frac{mn^2}{\sqrt{W}} + mn)$
which attains the desired lower bound.

Next we consider LAPACK's DGEQRF.
In the worst case, a panel factorization by DGEQR2 will incur
one slow memory access per arithmetic operation,
and so $O(\frac{n}{b} \cdot mb^2 ) =  O(mnb)$ for all panel factorizations.
For the overall algorithm to be guaranteed to attain
minimal bandwidth, we need $mnb = O(\frac{mn^2}{\sqrt{W}})$,
or $b = O(\frac{n}{\sqrt{W}})$. Since $b$ must also be at 
least about $\sqrt{W}$, this means $W = O(n)$,
or that fast memory size may be at most large enough
to hold a few rows of the matrix, or may be much smaller.

RGEQR3 does not alway minimize latency. For example,
considering applying RGEQR3 to a single panel
with $n=\sqrt{W}$ columns and $m>W$ rows, stored
in a block-column layout with $\sqrt{W}$-by-$\sqrt{W}$
blocks stored columnwise, as above. Then a recurrence
for the number of messages RGEQR3 requires is
\begin{eqnarray*}
\label{eqn:RGEQR3_latency}
L(m,n) & = & \left\{ \begin{array}{ll} 
             L(m,\frac{n}{2}) + L(m-\frac{n}{2},\frac{n}{2}) + 
             O(\frac{m}{\sqrt{W}})
             & {\rm if} \; n>1 \\
             O(\frac{m}{\sqrt{W}}) & {\rm if} \; n = 1  
             \end{array} \right. \nonumber \\
       & = & O(\frac{mn}{\sqrt{W}}) = O(m) \; {\rm when} \; n = \sqrt{W} 
\end{eqnarray*}
which is larger than the minimum $O(\frac{mn}{W}) = O(\frac{m}{\sqrt{W}})$
attained by sequential TSQR when $n = \sqrt{W}$.

In contrast to DGEQRF, RGEQRF, and RGEQR3, 
CAQR minimizes flops, bandwidth and latency
for all values of $W$.

%%%%%%%%%%%%%%%%%%%%%%%%%%%%%%%%%%%%%%%%%%%%%%%%%%%%%%%%%%%%%%%%%%%%%%
%%%%%%%%%%%%%%%%%%%%%%%%%%%%%%%%%%%%%%%%%%%%%%%%%%%%%%%%%%%%%%%%%%%%%%

\newcommand{\nunder}{{\underline{n}}}
\newcommand{\munder}{{\underline{m}}}

\section{Lower Bounds for TSQR}
\label{sec:LowerBounds_TSQR}

We present communication lower bounds for TSQR.
As we already mentioned for the sequential case, 
it is obviously necessary to read $mn$ words
from from slow to fast memory (the input),
and write $mn$ words from fast to slow memory (the output),
for a lower bound of $2mn$ words moved. Sequential
TSQR attains this trivial lower bound.
Since the size of a message is bounded by the size of
fast memory $W$, it clearly requires at least $\frac{2mn}{W}$
messages to send this much data. Since TSQR sends
$\frac{2mn}{\widetilde{W}} = \frac{2mn}{W - \frac{n(n+1)}{2}} \stackrel{<}{\approx} \frac{3mn}{W}$
messages, it attains this bound to within a constant factor, and is very close when
$W \gg n^2$.

For parallel TSQR, the lower bound on latency is obviously $\log P$, 
since TSQR needs to compute a nontrivial function of data that
is spread over $P$ processors, and a binary reduction tree
of depth $\log P$ clearly minimizes latency (by using the
butterfly variant). Parallel TSQR attains this lower bound too.
 
Bandwidth lower bounds for parallel TSQR are more interesting.
We analyze this in a way that applies to more general situations,
starting with the following:
Suppose processor 1 and processor 2 each own some of the arguments
of a function $f$ that processor 1 wants to compute. What is the least
volume of communication required to compute the function?
We are interested in smooth functions of real or complex arguments,
and so will use techniques from calculus rather than modeling
the arguments as bit strings.

In this way, we will derive necessary conditions on the function $f$
for it to be evaluable by communicating fewer than all of its arguments
to one processor. We will apply these conditions to various linear
algebra operations to capture our intuition that it is in fact necessary 
to move all the arguments to one processor for correct evaluation of $f$:
Subsection~\ref{sec:LowerBounds_TSQR_ss1} will show that
if $f$ is a bijection as a function of the $n$ arguments on processor 2,
and if processor 2 can only send one message to processor 1, then it 
indeed has to send all $n$ arguments
(part 3 of Lemma~\ref{lemma:bijection}).
Subsection~\ref{sec:LowerBounds_TSQR_ss2} extends this to reduction
operations where each processors sends one message to its parent in
a reduction tree, which is the case we are considering in this paper.
Subsection~\ref{sec:LowerBounds_TSQR_ss3} goes a step further and
asks whether less data can be sent overall by allowing processors 1 and 2
to exchange multiple but smaller messages; the answer is sometimes yes, but
again not for the reduction operations we consider.

\subsection{Communication lower bounds for one-way communication between
2 processors}
\label{sec:LowerBounds_TSQR_ss1}

Suppose $x^{(m)} \in \RR^m$ is owned by processor 1 (P1) and
$y^{(n)} \in \RR^n$ is owned by P2; we use superscripts to
remind the reader of the dimension of each vector-valued variable or function.
Suppose P1 wants to compute  
$f^{(r)}(x^{(m)},y^{(n)}): \RR^{m} \times \RR^{n} \rightarrow \RR^r$.
We first ask how much information P2 has to send to P1, assuming
it is allowed to send one message, consisting of $\nunder \leq n$ real
numbers, which themselves could be functions of $y^{(n)}$.
In other words, we ask if functions $h^{(\nunder )} (y^{(n)}): \RR^{n} \rightarrow \RR^{\nunder}$
and $F^{(r)} (x^{(m)}, z^{(\nunder)})  : \RR^{m} \times \RR^{\nunder} \rightarrow \RR^r$, 
exist such that 
$f^{(r)}(x^{(m)},y^{(n)}) = F^{(r)} (x^{(m)}, h^{(\nunder)} (y^{(n)}))$.
When $\nunder = n$, the obvious choice is to send the original data $y^{(n)}$,
so that $h^{(\nunder)} (y^{(n)}) = y^{(n)}$ is the identity function and
$f^{(r)} = F^{(r)}$. The interesting question is whether we can send
less information, i.e. $\nunder < n$.

Unless we make further restrictions on the function $h$ we are allowed
to use, it is easy to see that we can always choose $\nunder =1$, i.e. send
the least possible amount of information: We do this by using a 
space-filling curve \cite{sagan1994space} to represent each $y^{(n)} \in \RR^{(n)}$ by
one of several preimages $\tilde{y} \in \RR$. In other words, 
$h^{(1)} (y^{(n)})$ 
maps $y^{(n)}$ to a scalar $\tilde{y}$ that P1 can map back to 
$y^{(n)}$ by a space filling curve. 
This is obviously unreasonable, since it implies we could try to 
losslessly compress $n$ 64-bit floating point numbers into one 64-bit 
floating point number. 
However, by placing some reasonable smoothness restrictions on  
the functions we use, since we can only hope to evaluate (piecewise) smooth
functions in a practical way anyway, we will see that we can draw useful 
conclusions about practical computations. 
To state our results, we use the notation $J_x f(x,y)$ to denote the
$r \times m$ Jacobian matrix of $f^{(r)}$ with respect to the arguments
$x^{(m)}$.  Using the above notation, we state

\lemma{
Suppose it is possible to compute 
$f^{(r)}(x^{(m)},y^{(n)})$ on P1 by communicating $\nunder < n$ words
$h^{(\nunder )} (y^{(n)})$ from P2 to P1, and evaluating
$f^{(r)}(x^{(m)},y^{(n)}) = F^{(r)} (x^{(m)}, h^{(\nunder )} (y^{(n)}))$.
Suppose $h^{(\nunder )}$ and $F^{(r)}$ are continuously differentiable
on open sets. Then necessary conditions for this to be possible
are as follows.
\begin{enumerate}
\item Given any fixed $y^{(n)}$ in the open set, then for all
$x^{(m)}$ in the open set, 
the rows of $J_y f(x,y)$ must lie
in a fixed subspace of $\RR^n$
of dimension at most $\nunder < n$.
\item Given any fixed $\tilde{y}^{(\nunder)} \in \RR^{\nunder}$ satisfying
$\tilde{y}^{(\nunder)} = h^{(\nunder)} (y^{(n)})$ for some $y^{(n)}$ in
the interior of the open set, there is 
a set $C \subset \RR^{n}$ containing $y^{(n)}$,
of dimension at least $n-\nunder$,
such that for each 
$x$, $f(x,y)$ is constant for $y \in C$.
\item If $r=n$, and for each fixed $x$, $f^{(r)}(x,y^{(n)})$ is a bijection,
then it is necessary and sufficient to send $n$ words from P2 to P1 
to evaluate $f$.
\end{enumerate}
}
\label{lemma:bijection}

\rm 
\begin{proof}
Part 1 is proved simply by differentiating, using the chain rule,
and noting the dimensions of the Jacobians being multiplied:
\[
J_y^{(r \times n)} f^{(r)}(x,y) = J_h^{(r \times \nunder)} F^{(r)}(x,h) 
\cdot J_y^{(\nunder \times n)} h^{(\nunder)} (y)
\]
implying that for all $x$, each row of $J_y^{(r \times m)} f^{(r)} (x,y)$
lies in the space spanned by the $\nunder$ rows of 
$J_y^{(\nunder \times n)} h^{(\nunder)} (y)$.

Part 2 is a consequence of the implicit function theorem.
Part 3 follows from part 2, since if the function is a bijection,
then there is no set $C$ along which $f$ is constant.
\end{proof} \rm

Either part of the lemma can be used to derive 
lower bounds on the volume of communication needed to compute $f(x,y)$,
for example
by choosing an $\nunder$ equal to the lower bound minus 1, and 
confirming that either necessary condition in the Lemma is
violated, at least in some open set. 

We illustrate this for a simple matrix factorization problem.

\corollary{
Suppose P1 owns the $r_1 \times c$ matrix $A_1$, and
P2 owns the $r_2 \times c$ matrix $A_2$, with
$r_2 \geq c$. Suppose P1 wants to compute the
$c \times c$ Cholesky factor $R$ of
$R^T \cdot R = A_1^T \cdot A_1 + A_2^T \cdot A_2$,
or equivalently the $R$ factor in the $QR$ decomposition
of $\bmat{c} A_1 \\ A_2 \emat$. Then P2 has to communicate
at least $c(c+1)/2$ words to P1, and it is possible to
communicate this few, namely either the entries
on and above the diagonal of the symmetric $c \times c$ matrix $A_2^T \cdot A_2$, 
or the entries of its
Cholesky factor $R$, so that $R^T \cdot R = A_2^T \cdot A_2$ 
(equivalently, the $R$ factor of the $QR$ factorization of $A_2$).
}

\rm
\begin{proof}
That it is sufficient to communicate the $c(c+1)/2$ entries described
above is evident. We use  Corollary~1 to prove that these many words are
necessary. We use the fact that mapping between 
the entries on and above the diagonal of the symmetric positive definite
matrix and its Cholesky factor is a bijection
(assuming positive diagonal entries of the Cholesky factor).
To see that for any fixed $A_1$, $f(A_1,R) = $ the Cholesky factor
of $A_1^T \cdot A_1 + R^T \cdot R$ is a bijection,
note that it is a composition of three bijections:
the mapping from $R$ to the entries on and above the 
diagonal of $Y = A_2^T \cdot A_2$, the entries on and
above the diagonal of $Y$ and those on and above the diagonal
of $X = A_1^T \cdot A_1 + Y$, and the mapping between the entries on
and above the diagonal of $X$ and its Cholesky factor $f(A_1,R)$.
\end{proof} \rm

\subsection{Reduction operations}
\label{sec:LowerBounds_TSQR_ss2}

We can extend this result slightly to make it apply to the case of
more general reduction operations, where one processor P1 is trying to
compute a function of data initially stored on multiple other processors
P2 through P$s$. We suppose that there is a tree
of messages leading from these processors eventually reaching P1.
Suppose each P$i$ only sends data up the tree, so that the communication
pattern forms a DAG (directed acylic graph) with all paths ending at P1. 
Let P$i$'s data be denoted $y^{(n)}$.
Let all the variables on P1 be denoted $x^{(m)}$,
and treat all the other variables on the other processors as constants.
Then exactly the same analysis as above applies, and we can conclude that
{\em every} message along the unique path from P$i$ to P1 has the same 
lower bound on its size, as determined by Lemma~1.
This means Corollary~1 extends to include reduction operations where
each operation is a bijection between one input (the other being fixed)
and the output. In particular, it applies to TSQR.

We emphasize again that using a real number model to draw conclusions about
finite precision computations must be done with care. For example,
a bijective function depending on many variables could hypothetically 
round to the same floating point output for all floating point inputs,
eliminating the need for any communication or computation 
for its evaluation. But this is not the
case for the functions we are interested in.

Finally, we note that the counting must be done slightly 
differently for the QR decomposition of complex data,
because the diagonal entries $R_{i,i}$ are generally
taken to be real. Alternatively, there is a degree of
freedom in choosing each row of $R$, which can be
multiplied by an arbitrary complex number of absolute
value 1.

\subsection{Extensions to two-way communication}
\label{sec:LowerBounds_TSQR_ss3}

While the result of the previous subsection is adequate for the results 
of this paper,
we note that it may be extended as follows. For motivation, suppose that
P1 owns the scalar $x$, and wants to evaluate the polynomial
$\sum_{i=1}^{n} y_i x^{i-1}$, where P2 owns the vector $y^{(n)}$.
The above results can be used to show that P2 needs to send $n$
words to P1 (all the coefficients of the polynomial, for example).
But there is an obvious way to communicate just 2 words:
(1) P1 sends $x$ to P2, (2) P2 evaluates the polynomial, and
(3) P2 sends the value of the polynomial back to P1.

More generally, one can imagine $k$ phases, during each of which
P1 sends one message to P2 and then P2 sends one message to P1.
The contents of each message can be any smooth functions of all 
the data available to the sending processor, either originally 
or from prior messages. At the end of the $k$-th phase, P1 then
computes $f(x,y)$.

More specifically, the computation and communication proceeds as
follows:
\begin{itemize}
\item In Phase 1, P1 sends $g_1^{(m_1)} (x^{(m)})$ to P2
\item In Phase 1, P2 sends $h_1^{(n_1)} (y^{(n)}, g_1^{(m_1)} (x^{(m)}))$ to P1
\item In Phase 2, P1 sends $g_2^{(m_2)} (x^{(m)}, h_1^{(n_1)} (y^{(n)}, g_1^{(m_1)} (x^{(m)})))$ to P2
\item In Phase 2, P2 sends 
$h_2^{(n_2)} (y^{(n)}, 
g_1^{(m_1)} (x^{(m)}),
g_2^{(m_2)} (x^{(m)}, h_1^{(n_1)} ( y^{(n)}, g_1^{(m_1)} (x^{(m)}))))$ to P1
\item $\dots$
\item In Phase $k$, P1 sends
$g_k^{(m_k)} (x^{(m)}, h_1^{(n_1)} (\dots), h_2^{(n_2)} (\dots) , \dots , 
h_{k-1}^{(n_{k-1})} (\dots) )$ to P2
\item In Phase $k$, P2 sends
$h_k^{(n_k)} (y^{(n)}, g_1^{(m_1)} (\dots), g_2^{(m_2)} (\dots) , \dots , 
g_{k}^{(m_k)} (\dots) )$ to P1
\item P1 computes 
\begin{eqnarray*}
f^{(r)} (x^{(m)} , y^{(n)} ) & = &
F^{(r)} ( x^{(m)}, h_1^{(n_1)} ( y^{(n)}, g_1^{(m_1)} (x^{(m)})), 
\\ & & 
h_2^{(n_2)} ( y^{(n)}, 
g_1^{(m_1)} (x^{(m)}), 
g_2^{(m_2)} (x^{(m)}, h_1^{(n_1)} ( y^{(n)}, g_1^{(m_1)} (x^{(m)})))),
\\ & & \dots
\\ & & h_k^{(n_k)} (y^{(n)}, g_1^{(m_1)} (\dots), g_2^{(m_2)} (\dots) , \dots , 
g_{k}^{(m_k)} (\dots) ))
\end{eqnarray*}
\end{itemize}

\lemma{
Suppose it is possible to compute 
$f^{(r)}(x^{(m)},y^{(n)})$ on P1 by the scheme described above.
Suppose all the functions involved are continuously differentiable
on open sets. Let $\nunder = \sum_{i=1}^k n_i$ and
$\munder = \sum_{i=1}^k m_i$.
Then necessary conditions for this to be possible
are as follows.
\begin{enumerate}
\item 
Suppose $\nunder < n$ and $\munder \leq m$, ie. P2 cannot communicate
all its information to P1, but P1 can potentially send its information
to P2.
Then there is a set $C_x \subset \RR^m$ 
of dimension at least $m-\munder$ 
and a set $C_y \subset \RR^n$ 
of dimension at least $n-\nunder$ such that
for $(x,y) \in C = C_x \times C_y$, 
the value of $f(x,y)$ is independent of $y$.
\item 
If $r=n=m$, and for each fixed $x$ or fixed $y$, 
$f^{(r)}(x^{(m)},y^{(n)})$ is a bijection,
then it is necessary and sufficient to send $n$ words from P2 to P1 
to evaluate $f$.
\end{enumerate}
}
\rm

\begin{proof}
We define the sets $C_x$ and $C_y$ by the following constraint equations,
one for each communication step in the algorithm:
\begin{itemize}
\item 
$\tilde{g}_1^{(m_1)} = g_1^{(m_1)} (x^{(m)})$ is a fixed constant, 
placing $m_1$ smooth constraints on $x^{(m)}$.
\item 
In addition to the previous constraint,
$\tilde{h}_1^{(n_1)} = h_1^{(n_1)} (y^{(n)}$, $g_1^{(m_1)} (x^{(m)}))$ 
is a fixed constant,
placing $n_1$ smooth constraints on $y^{(n)}$.
\item 
In addition to the previous constraints, \linebreak
$\tilde{g}_2^{(m_2)} = g_2^{(m_2)} (x^{(m)}, h_1^{(n_1)} (y^{(n)}, g_1^{(m_1)} (x^{(m)})))$ 
is a fixed constant, placing $m_2$ more smooth constraints on $x^{(m)}$.
\item 
In addition to the previous constraints, \linebreak
$\tilde{h}_2^{(n_2)} = h_2^{(n_2)} (y^{(n)}, 
g_1^{(m_1)} (x^{(m)}),
g_2^{(m_2)} (x^{(m)}, h_1^{(n_1)} ( y^{(n)}, g_1^{(m_1)} (x^{(m)}))))$ 
is a fixed constant,
placing $n_2$ more smooth constraints on $y^{(n)}$.
\item 
\dots
\item 
In addition to the previous constraints, \linebreak
$\tilde{g}_k^{(m_k)} = g_k^{(m_k)} (x^{(m)}, h_1^{(n_1)} (\dots), h_2^{(n_2)} (\dots) , \dots , 
h_{k-1}^{(n_{k-1})} (\dots) )$ is a fixed constant,
placing $m_k$ more smooth constraints on $x^{(m)}$.
\item 
In addition to the previous constraints, \linebreak
$\tilde{h}_k^{(n_k)} = h_k^{(n_k)} (y^{(n)}, g_1^{(m_1)} (\dots), g_2^{(m_2)} (\dots) , \dots , 
g_{k}^{(m_k)} (\dots) )$ is a fixed constant,
placing $n_k$ more smooth constraints on $y^{(n)}$.
\end{itemize}
Altogether, we have placed 
$\nunder = \sum_{i=1}^k n_i < n$ smooth constraints on $y^{(n)}$ and
$\munder = \sum_{i=1}^k m_i \leq m$ smooth constraints on $x^{(m)}$,
which by the implicit function theorem define surfaces 
$C_y (\tilde{h}_1^{(n_1)} , \dots , \tilde{h}_k^{(n_k)} )$
and 
$C_x (\tilde{g}_1^{(m_1)} , \dots , \tilde{g}_k^{(m_k)} )$,
of dimensions at least
$n - \nunder > 0$ and $m - \munder \geq 0$, respectively,
and parameterized by
$\{\tilde{h}_1^{(n_1)} , \dots , \tilde{h}_k^{(n_k)} \}$ and
$\{\tilde{g}_1^{(m_1)} , \dots , \tilde{g}_k^{(m_k)} \}$,
respectively.
For $x \in C_x$ and $y \in C_y$, the values communicated by 
P1 and P2 are therefore constant. Therefore, for $x \in C_x$
and $y \in C_y$, $f(x,y) = F(x,h_1, \dots, h_k)$ depends only on $x$,
not on $y$. This completes the first part of the proof.

For the second part, we know that if $f(x,y)$ is a bijection
in $y$ for each fixed $x$, then by the first part
we cannot have $\nunder < n$, because otherwise 
$f(x,y)$ does not depend on $y$ for certain values of $x$,
violating bijectivity. But if we can send $\nunder = n$ words from
P2 to P1, then it is clearly possible to compute $f(x,y)$ by
simply sending every component of $y^{(n)}$ from P2 to P1 explicitly.
\end{proof}

\rm

\corollary{
Suppose P1 owns the $c$-by-$c$ upper triangular matrix $R_1$,
and P2 owns the $c$-by-$c$ upper triangular matrix $R_2$, and
P1 wants to compute the R factor in the QR decomposition of
$\bmat{c} R_1 \\ R_2 \emat$. Then it is necessary and sufficient
to communicate $c(c+1)/2$ words from P2 to P1 (in particular, the
entries of $R_2$ are sufficient).
}

\rm 

We leave extensions to general communication patterns among
multiple processors to the reader.

% \end{document}

\section{Lower Bounds for CAQR}
\label{sec:LowerBounds_CAQR}

In this section, we review known lower bounds on communication 
bandwidth for parallel and sequential $\Theta (n^3)$ 
matrix-matrix multiplication 
of matrices stored in 
% 1-D and 
2-D layouts, 
extend some of them to the rectangular case, and then extend them to LU
and QR, showing that our sequential and parallel CAQR 
algorithms have optimal communication complexity with respect
to both bandwidth (in a Big-Oh sense, and sometimes modulo
polylogarithmic factors).

We will also use the simple fact that if $B$ is a lower bound on the number
of words that must be communicated to implement an algorithm,
and if $W$ is the size of the local memory (in the parallel case) 
or fast memory (in the sequential case),
so that $W$ is the largest possible size of a message,
then $B/W$ is a lower bound on the latency, i.e. the number
of messages needed to move $B$ words into or out of the memory.
We use this to derive lower bounds on latency, which are
also attained by our algorithms (again in a Big-Oh sense,
and sometimes modulo polylogarithmic factors).

We begin in section~\ref{SS:MMlowerbounds} by reviewing
known communication complexity bounds for 
$\Theta (n^3)$ matrix multiplication,
due first to Hong and Kung \cite{hong1981io} in the sequential case,
and later proved more simply and extended to the parallel case
by Irony, Toledo and Tiskin \cite{irony2004communication}.  

It is easy to extend lower bounds for matrix multiplication
to lower bounds for LU decomposition via the following 
reduction of matrix multiplication to LU:
\begin{equation}\label{eq:GEMM-to-LU}
\begin{pmatrix}
I & 0 & -B \\
A & I & 0 \\
0 & 0 & I \\
\end{pmatrix}
=
\begin{pmatrix}
I &   &   \\
A & I &   \\
0 & 0 & I \\
\end{pmatrix}
\begin{pmatrix}
I & 0 & -B \\
  & I & A \cdot B \\
  &   & I \\
\end{pmatrix}.
\end{equation}
See \cite{grigori2008calu} for an implementation of 
parallel LU that attains these bounds.
See \cite{toledo1997locality} for an implementation of sequential
LU and a proof that it attains the bandwidth lower bound
(whether the latency lower bound is attained is an open problem).

It is reasonable to expect that lower bounds for matrix multiplication
will also apply (at least in a Big-Oh sense) to other 
one-sided factorizations, such as QR.
As we will see, QR is not as simple as LU.

All this assumes commutative and associative
reorderings of conventional $\Theta(n^3)$ 
matrix multiplication,
and so excludes faster algorithms using distributivity
or special constants, such as those of
Strassen \cite{strassen1969gaussian} or Coppersmith and Winograd \cite{coppersmith1982asymptotic},
and their use in asymptotically fast versions of
LU and QR \cite{FastLinearAlgebraIsStable}.
Extending communication lower bounds to these asymptotically
faster algorithms is an open problem.

\subsection{Matrix Multiplication Lower Bounds}
\label{SS:MMlowerbounds}

We review lower bounds in \cite{hong1981io,irony2004communication} 
for multiplication of two $n$-by-$n$ matrices $C = A \cdot B$
using commutative and associative (but not distributive)
reorderings of the usual $\Theta(n^3)$ algorithm.
In the sequential case, they assume that $A$ and $B$ initially reside
in slow memory, that there is a fast memory of size $W < n^2$,
and that the product $C = A \cdot B$ must be computed and eventually 
reside in slow memory.
They bound from below the number of
words that need to be moved between slow memory and fast memory
to perform this task:
\begin{equation}\label{eqn_MatMul_seq_bw_lowerbound}
{\rm \#\ words\ moved} \geq \frac{n^3}{2 \sqrt{2} W^{1/2}} - W \approx
\frac{n^3}{2 \sqrt{2} W^{1/2}}  \; \; .
\end{equation}
Since only $W$ words can be moved in one message, this also provides
a lower bound on the number of messages:
\begin{equation}\label{eqn_MatMul_seq_lat_lowerbound}
{\rm \#\ messages} \geq \frac{n^3}{2 \sqrt{2} W^{3/2}} - 1 \approx
\frac{n^3}{2 \sqrt{2} W^{3/2}} \; \; .
\end{equation}
In the rectangular case, where $A$ is $n$-by-$r$, $B$ is $r$-by-$m$,
and $C$ is $n$-by-$m$, so that the number of arithmetic operations in
the standard algorithm is $2mnr$, the above two results still apply,
but with $n^3$ replaced by $mnr$.

The parallel case is considered in \cite{irony2004communication}.
There is
actually a spectrum of algorithms, from the so-called 2D case,
that use little extra memory beyond that needed to store equal fractions
of the matrices $A$, $B$ and $C$
(and so about $3n^2/P$ words for each of $P$ processors, in the square case), 
to the 3D case, where each input matrix is replicated up to $P^{1/3}$ times,
so with each processor needing memory of size $n^2/P^{2/3}$ in the square case.
We only consider the 2D case, which is the conventional, memory 
scalable approach.
In the 2D case, with square matrices, Irony et al show that
if each processor has $\mu n^2 /P$ words of local
memory, and $P \geq 32 \mu^3$, then at least one of the processors
must send or receive at least the following number of words:
\begin{equation}\label{eqn_MatMul_par_bw_lowerbound}
{\rm \#\ words\ sent\ or\ received} \geq \frac{n^2}{4 \sqrt{2} (\mu P)^{1/2}} 
\end{equation}
and so using at least the following number of messages
(assuming a maximum message size of $n^2/P$):
\begin{equation}\label{eqn_MatMul_par_lat_lowerbound}
{\rm \#\ messages} \geq \frac{P^{1/2}}{4 \sqrt{2} (\mu)^{3/2}} \; \; .
\end{equation}

We wish to extend this to the case of rectangular matrices.
We do this in preparation for analyzing CAQR in the rectangular case.
The proof is a simple extension of Thm.~4.1 in 
\cite{irony2004communication}.  

\theorem{
Consider the conventional matrix multiplication algorithm applied to
$C = A \cdot B$ where
$A$ is $n$-by-$r$, $B$ is $r$-by-$m$, and $C$ is $n$-by-$m$,
implemented on a $P$ processor distributed memory parallel computer.
Let $\bar{n}$, $\bar{m}$ and $\bar{r}$ be the sorted values of
$n$, $m$, and $r$, i.e. $\bar{n} \geq \bar{m} \geq \bar{r}$.
Suppose each processor has $3\bar{n}\bar{m}/P$ words of local memory,
so that it can fit 3 times as much as $1/P$-th of the largest of the
three matrices. Then as long as
\begin{equation}\label{eqn:RectMM}
\bar{r} \geq \sqrt{\frac{864 \bar{n}\bar{m}}{P}}
\end{equation}
(i.e. none of the matrices is ``too rectangular'')
then the number of words at least one processor must send or
receive is 
\begin{equation}\label{eqn:RectMM_bw}
{\rm \#\ words\ moved} \geq 
\frac{\sqrt{\bar{n}\bar{m}} \cdot \bar{r}}
{\sqrt{96 P}}
\end{equation}
and the number of messages is 
\begin{equation}\label{eqn:RectMM_lat}
{\rm \#\ messages} \geq 
\frac{ \sqrt{P} \cdot \bar{r}}
{ \sqrt{864 \bar{n}\bar{m}} } 
\end{equation}
}

\rm

\begin{proof}
We use (\ref{eqn_MatMul_seq_bw_lowerbound}) 
with $\bar{m} \bar{n} \bar{r}/P$ substituted for $n^3$,
since at least one processor does this much arithmetic,
and $W = 3\bar{n}\bar{m}/P$ words of local memory.
The constants in inequality (\ref{eqn:RectMM})
are chosen so that the first term in
(\ref{eqn_MatMul_seq_bw_lowerbound}) is at least $2W$,
and half the first term is a lower bound.
\end{proof}

It is well-known that the communication lower bound for 
sequential matrix multiplication is attained by
``tiling'' or ``blocking'' the matrices into square
blocks of dimension $\sqrt{W/3}$, and for parallel
matrix multiplication by Cannon's algorithm \cite{cannon1969cellular}.

\subsection{Lower Bounds for CAQR}\label{SS:lowerbounds:2d}

Now we need to extend our analysis of matrix multiplication.
We assume all variables are real; extensions to the complex
case are straightforward.
Suppose $A = QR$ is $m$-by-$n$, $n$ even, 
so that 
\[
\bar{Q}^T \cdot \bar{A} 
\equiv 
\left(
    Q( 1:m, 1:\frac{n}{2} )
\right)^T 
\cdot 
A( 1:m, \frac{n}{2}+1:n )
= 
R( 1:\frac{n}{2}, \frac{n}{2}+1:n )
\equiv 
\bar{R} \; \; .
\]
It is easy to see that $\bar{Q}$
depends only on the first $\frac{n}{2}$ columns of $A$, and
so is independent of $\bar{A}$. The obstacle to directly
applying existing lower bounds for matrix multiplication 
of course is that $\bar{Q}$ is not represented as an explicit 
matrix, and $\bar{Q}^T \cdot \bar{A}$ is not implemented by 
straightforward matrix multiplication. 
Nevertheless, we argue that the same data dependencies
as in matrix multiplication can be found inside many
implementations
of $\bar{Q}^T \cdot \bar{A}$,
and that therefore
the geometric ideas underlying the analysis in
\cite{irony2004communication} still apply.
Namely, there are two data structures $\tilde{Q}$
and $\tilde{A}$ indexed with pairs of subscripts
$(j,i)$ and $(j,k)$ respectively with the following properties.
\begin{itemize}
\item $\tilde{A}$ stores $\bar{A}$ as well as all intermediate
results which may overwrite $\bar{A}$. 
\item $\tilde{Q}$ represents $\bar{Q}$, i.e., an $m$-by-$\frac{n}{2}$
orthogonal matrix. Such a matrix is a member of the Stiefel manifold
of orthogonal matrices, and is known to require 
$\frac{mn}{2} - \frac{n}{4}(\frac{n}{2}+1)$ independent parameters
to represent, with column $i$ requiring $m-i$ parameters,
although a particular algorithm may represent
$\bar{Q}$ using more data.
\item The algorithm operates mathematically independently on each column
of $\bar{A}$, i.e., methods like that of Strassen are excluded.
This means that the algorithm performs at least 
$\frac{mn}{2} - \frac{n}{4}(\frac{n}{2}+1)$ multiplications 
on each $m$-dimensional column vector of $\bar{A}$
(see subsection~\ref{SS:lowerbounds:2d:flops} for a proof),
% with at least
% $m-i$ multiplications to get $\bar{R}_{i,k}$ 
and does the same operations on each column of $\bar{A}$. 
\item For each $(i,k)$ indexing $\bar{R}_{i,k}$, which is the 
component of the $k$-th column $\bar{A}_{:,k}$ of $\bar{A}$ in
the direction of the $i$-th column $\bar{Q}_{:,i}$ of $\bar{Q}$,
it is possible to identify at least $m-i$ common components of $\tilde{A}_{:,k}$
and of $\tilde{Q}_{:,i}$ such that a parameter associated with
$\tilde{Q}_{j,i}$ is multiplied by a value stored in $\tilde{A}_{j,k}$.
\end{itemize}
The last point, which says that $\bar{Q}^T \cdot \bar{A}$
has at least the same dependencies as matrix multiplication, 
requires illustration.
\begin{itemize}
\item Suppose $\bar{Q}$ is represented as a product of $\frac{n}{2}$
Householder reflections with a projection $\hat{Q}$ onto the
first $\frac{n}{2}$ coordinates,
$\bar{Q} = 
(I - \tau_1 u_1u_1^T)
	\cdots 
(I - \tau_{n/2} u_{n/2} u_{n/2}^T) 
\hat{Q}$,
normalized in the 
conventional way where the topmost nonzero entry of each $u_j$ is one,
and $\hat{Q}$ consists of the first $n/2$ columns of the $n$-by-$n$
identity matrix.
Then $\tilde{Q}_{j,i} = u_i(j)$ is multiplied by some intermediate value of
$\bar{A}_{j,k}$, i.e. $\tilde{A}_{j,k}$.
\item Suppose $\bar{Q}$ is represented as a product of block
Householder transformations $(I-Z_1U_1^T) \cdots (I-Z_f U_f^T) \hat{Q}$
where $U_g$ and $Z_g$ are $m$-by-$b_g$ matrices, $U_g$ consisting
of $b_g$ Householder vectors side-by-side. 
Again associate $\tilde{Q}_{j,i}$ with
the $j$-th entry of the $i$-th Householder vector $u_i(j)$.
\item Recursive versions of QR \cite{elmroth1998new} apply 
blocked Householder transformations organized so as to better
use BLAS3, but still let us use the approach of the last bullet.
\item Suppose $\bar{Q}$ is represented as a product of 
$\frac{mn}{2} - \frac{n}{4}(\frac{n}{2}+1)$ Givens rotations,
each one creating a unique subdiagonal zero entry in $A$ which is
never filled in.  There are
many orders in which these zeros can be created, and possibly
many choices of row that each Givens rotation may rotate with to zero
out its desired entry. 
If the desired zero entry in $A_{j,i}$ is created by
the rotation in rows $j'$ and $j$, $j'<j$, 
then associate $\tilde{Q}_{j,i}$ with the value
of the cosine in the Givens rotation, since this will be multiplied
by $\bar{A}_{j,k}$.
\item Suppose, finally, that we use CAQR to perform the
QR decomposition, so that $\bar{Q} = Q_1 \cdots Q_f \hat{Q}$, where
each $Q_g$ is the result of TSQR on $b_g$ columns.
Consider without loss of generality $Q_1$, which operates
on the first $b_1$ columns of $A$.
We argue that TSQR still produces $m-i$ parameters associated
with column $i$ as the above methods. Suppose there are $P$
row blocks, each of dimension $\frac{m}{P}$-by-$b_1$.
Parallel TSQR initially does QR independently on each block, using
any of the above methods; we associate multipliers as 
above with the subdiagonal entries in each block. Now
consider the reduction tree that combines $q$ different $b_1$-by-$b_1$
triangular blocks at any particular node.
%As described in 
% subsection~\ref{SS:TSQR:localQR:structured}, 
%subsection~\ref{sec:TSQR_optimal_localopt}
%this generates
This generates
$(q-1)b_1(b_1+1)/2$ parameters that multiply the equal number of entries 
of the $q-1$ triangles being zeroed out, and so can be associated with
appropriate entries of $\tilde{Q}$. Following the reduction tree, we
see that parallel TSQR produces exactly 
as many parameters as Householder reduction,
and that these may be associated one-for-one with all subdiagonal
entries of $\tilde{Q}(:,1:b_1)$ and $\tilde{A}(:,1:b_1)$ as above.
Sequential TSQR reduction is analogous.
\end{itemize}
We see that we have only tried to capture the dependencies
of a fraction of the arithmetic operations performed by various
QR implementations; this is all we need for a lower bound.

Now we resort to the geometric approach of 
\cite{irony2004communication}: Consider a three dimensional 
block of lattice points, indexed by $(i,j,k)$.
Each point on the $(i,0,k)$ face is associated with $\bar{R}_{i,k}$,
for $1 \leq i,k \leq \frac{n}{2}$.
Each point on the $(0,j,k)$ face is associated with $\tilde{A}_{j,k}$,
for $1 \leq k \leq \frac{n}{2}$ and $1 \leq j \leq m$.
Each point on the $(i,j,0)$ face is associated with $\tilde{Q}_{j,i}$,
for $1 \leq i \leq \frac{n}{2}$ and $1 \leq j \leq m$.
Finally, each interior point $(i,j,k)$ for 
$1 \leq i,k \leq \frac{n}{2}$ and $1 \leq j \leq m$ represents the
multiplication $\tilde{Q}_{j,i} \cdot \tilde{A}_{j,k}$.
The point is that the multiplication at $(i,j,k)$ cannot occur
unless $\tilde{Q}_{j,i}$ and $\tilde{A}_{j,k}$ are together in memory.

Finally, we need the Loomis-Whitney inequality \cite{loomis1949inequality}:
Suppose $V$ is a set of lattice points in 3D, 
$V_i$ is projection of $V$ along $i$ onto the $(j,k)$ plane,
and similarly for $V_j$ and $V_k$. Let $|V|$ denote the cardinality of $V$,
i.e. counting lattice points. Then 
$|V|^2 \leq |V_i| \cdot |V_j| \cdot |V_k|$. 
We can now state

\lemma{
Suppose a processor with local (fast) memory of size $W$ is participating
in the QR decomposition of an $m$-by-$n$ matrix, $m \geq n$, using an
algorithm of the sort discussed above.
There may or may not be other processors participating (i.e. this lemma covers
the sequential and parallel cases). Suppose the processor performs $F$
multiplications. Then the processor must move
the following number of words into or out of its memory:
\begin{equation}\label{Thm:1_bw}
{\rm \#\ of\ words\ moved} \geq \frac{F}{(8W)^{1/2}} - W
\end{equation}
using at least the following number of messages:
\begin{equation}\label{Thm:1_lat}
{\rm \#\ of\ messages} \geq \frac{F}{(8W^3)^{1/2}} - 1
\end{equation}
}
\label{lemma:LB}

\rm
\begin{proof}
The proof closely follows that of Lemma~3.1 in \cite{irony2004communication}.
We decompose the computation into phases. Phase $l$ begins when the
total number of words moved into and out of memory is exactly $lW$.
Thus in each phase, except perhaps the last, the memory loads and stores
exactly $W$ words.

The number of words $n_A$ from different $\tilde{A}_{jk}$ that the processor
can access in its memory during a phase is $2W$, since each word was either
present at the beginning of the phase or read during the phase.
Similarly the number of coefficients $n_Q$ from different $\tilde{Q}_{ji}$
also satisfies $n_Q \leq 2W$. Similarly, the number $n_R$ of locations
into which intermediate results like $\tilde{Q}_{ji} \cdot \tilde{A}_{jk}$
can be accumulated or stored is at most $2W$. Note that these intermediate
results could conceivably be stored or accumulated in $\tilde{A}$ because 
of overwriting; this does not affect the upper bound on $n_R$.

By the Loomis-Whitney inequality, the maximum number of useful multiplications
that can be done during a phase (i.e. assuming intermediate results are
not just thrown away) is bounded by 
$\sqrt{n_A \cdot n_Q \cdot n_R} \leq \sqrt{8W^3}$. Since the processor does
$F$ multiplications, the number of full phases required is at least
\[
\left\lfloor \frac{F}{\sqrt{8W^3}} \right\rfloor \geq \frac{F}{\sqrt{8W^3}} -1
\]
so the total number of words moved is $W$ times larger, i.e. at least
\[
{\rm \#\ number\ of\ words\ moved} \geq
\frac{F}{\sqrt{8W}} -W \; \; .
\]
The number of messages follows by dividing by $W$, the maximum
message size.
\end{proof}
\rm

% \subsubsection{Sequential CAQR}
% \label{sec:SeqCAQR}

The following is our main result for sequential CAQR:
\corollary{
Consider a single processor computing the QR decomposition of
an $m$-by-$n$ matrix with $m \geq n$, using an algorithm of the
sort discussed above. Then the number of words moved between
fast and slow memory is at least
\begin{equation}\label{Thm:2_bw}
{\rm \#\ of\ words\ moved} \geq \frac
{\frac{mn^2}{4}-\frac{n^2}{8}(\frac{n}{2}+1)}
{(8W)^{1/2}} - W
\geq \frac{3n^2(m - \frac{4}{3})}{16(8W)^{1/2}} - W
\end{equation}
using at least the following number of messages:
\begin{equation}\label{Thm:2_lat}
{\rm \#\ of\ messages} \geq \frac
{\frac{mn^2}{4}-\frac{n^2}{8}(\frac{n}{2}+1)}
{(8W^3)^{1/2}} - 1
\geq \frac{3n^2(m - \frac{4}{3})}{16(8W^3)^{1/2}} - 1
\end{equation}
}
\label{corollary:SeqCAQR}
\rm
\begin{proof}
The proof follows easily from Lemma~\ref{lemma:LB} by
using the lower bound
$F \geq {\frac{mn^2}{4}-\frac{n^2}{8}(\frac{n}{2}+1)}$ on the number
of multiplications by any algorithm in the class discussed above
(see Lemma~\ref{lemma:F_lowerbound} in 
subsection~\ref{SS:lowerbounds:2d:flops} for a proof).
\end{proof}

\rm
The lower bound could be increased by a constant factor by using
a specific number of multiplications (say $mn^2 - n^3 / 3$
using Householder reductions), instead of arguing more generally based on
the number of parameters needed to represent orthogonal matrices.

Comparing to the performance model in
Section~\ref{sec:CAQR_sequential}, especially
Table~\ref{tbl:CAQR:seq:model:opt}, 
we see that sequential CAQR attains these bounds to within a 
constant factor.
% Comparing to Equation \eqref{eq:CAQR:seq:modeltime:P:opt} in Appendix
% \ref{S:CAQR-seq-detailed} or the presentation in Section \ref{S:CAQR-seq},
% we see that CAQR attains these bounds to within a constant factor.

% \subsubsection{Parallel CAQR}

The following is our main result for parallel CAQR:
\corollary{
Consider a parallel computer with $P$ processors
and $W$ words of memory per processor
computing the QR decomposition of
an $m$-by-$n$ matrix with $m \geq n$, using an algorithm of the
sort discussed above. 
Then the number of words sent and received
by at least one processor
is at least
\begin{equation}\label{Thm:3a_bw}
{\rm \#\ of\ words\ moved} \geq \frac
{\frac{mn^2}{4}-\frac{n^2}{8}(\frac{n}{2}+1)}
{P(8W)^{1/2}} - W
\geq \frac{3n^2(m - \frac{4}{3})}{16p(8W)^{1/2}} - W
\end{equation}
using at least the following number of messages:
\begin{equation}\label{Thm:3a_lat}
{\rm \#\ of\ messages} \geq \frac
{\frac{mn^2}{4}-\frac{n^2}{8}(\frac{n}{2}+1)}
{P(8W^3)^{1/2}} - 1
\geq \frac{3n^2(m - \frac{4}{3})}{16P(8W^3)^{1/2}} - 1
\end{equation}
In particular, when each processor has $W = mn/P$ words of memory
and the matrix is not too rectangular, $n \geq \frac{2^{11}m}{P}$,
then the number of words sent and received
by at least one processor is at least
\begin{equation}\label{Thm:3b_bw}
{\rm \#\ of\ words\ moved} \geq 
\sqrt{\frac{m n^3}{2^{11}P}}
\end{equation}
using at least the following number of messages:
\begin{equation}\label{Thm:3b_lat}
{\rm \#\ of\ messages} \geq 
\sqrt{\frac{nP}{2^{11}m}} \; \; .
\end{equation}
In particular, in the square case $m=n$, we get that
as long as $P \geq 2^{11}$, 
then the number of words sent and received
by at least one processor is at least
\begin{equation}\label{Thm:4b_bw}
{\rm \#\ of\ words\ moved} \geq 
{\frac{n^2}{2^{11/2}P^{1/2}}}
\end{equation}
using at least the following number of messages:
\begin{equation}\label{Thm:4b_lat}
{\rm \#\ of\ messages} \geq 
\sqrt{\frac{P}{2^{11}}} \; \; .
\end{equation}
}
\rm
\begin{proof}
The result follows from the previous Corollary, since
at least one processor has to do $1/P$-th of the work.
\end{proof}

Comparing to the performance model in Section~\ref{sec:CAQR_parallel},
especially Table~\ref{tbl:CAQR:par:model:opt}, we see that
parallel CAQR attains these bounds to within a constant factor.
% Comparing to Equations \eqref{eq:CAQR:par:opt:lat} and
% \eqref{eq:CAQR:par:opt:bw} in Section \ref{SS:CAQR:par:opt}, we see
% that CAQR attains these bounds to within a polylog factor.

\subsection{Lower Bounds on Flop Counts for QR}
\label{SS:lowerbounds:2d:flops}

This section proves lower bounds on arithmetic for {\em any} ``columnwise'' 
implementation of QR, by which we mean one whose operations can be reordered
so as to be left looking, i.e. the operations that compute columns $i$
of $Q$ and $R$ depend on data only in columns 1 through $i$ of $A$.
The mathematical dependencies are such that columns $i$ of $Q$ and $R$
do only depend on columns 1 through $i$ of $A$, but saying that operations
only depend on these columns eliminates algorithms like Strassen.
(It is known that QR can be done asymptotically as fast as any fast
matrix multiplication algorithm like Strassen, and stably 
\cite{FastLinearAlgebraIsStable}.)

This section says where the lower bound on
$F$ comes from that is used in the proof of 
Corollary~\ref{corollary:SeqCAQR} above.

The intuition is as follows. 
Suppose $A = QR$ is $m$-by-$(j+1)$, 
so that 
\[
\bar{Q}^T \cdot \bar{A} 
\equiv (Q(1:m,1:j))^T \cdot A(1:m,j+1) = 
R(1:j, j+1)
\equiv \bar{R} \; \; .
\]
where $\bar{Q}$ only depends on the first $j$ columns of $A$, and
is independent of $\bar{A}$. As an arbitrary $m$-by-$j$ orthogonal 
matrix, a member of the Stiefel manifold of dimension
$mj-j(j+1)/2$, $\bar{Q}$ requires $mj-j(j+1)/2$ independent
parameters to represent. We will argue that no matter how $\bar{Q}$
is represented, i.e. without appealing to the special structure of
Givens rotations or Householder transformations, that unless
$mj-j(j+1)/2$ multiplications are performed to compute $\bar{R}$
it cannot be computed correctly, because it cannot depend on
enough parameters. 

Assuming for a moment that this is true, we get a lower bound
on the number of multiplications needed for QR on an $m$-by-$n$ matrix 
by summing 
$\sum_{j=1}^{n-1} [mj-j(j+1)/2] = \frac{mn^2}{2} - \frac{n^3}{6} + O(mn)$.
The two leading terms are half the multiplication count for Householder
QR (and one fourth of the total operation count, including additions).
So the lower bound is rather tight.

Again assuming this is true, we get a lower bound on the
value $F$ in Corollary~\ref{corollary:SeqCAQR} by multiplying 
$\frac{n}{2} \cdot (m\frac{n}{2} - \frac{n}{2}(\frac{n}{2}+1)/2)
= \frac{mn^2}{4} - \frac{n^2}{8}(\frac{n}{2}+1) \leq F$.

Now we prove the main assertion, that $mj-j(j+1)/2$ multiplications are needed
to compute the single column $\bar{R} = \bar{Q}^T \cdot \bar{A}$, no matter how 
$\bar{Q}$ is represented. We model the computation as a DAG (directed
acyclic graph) of operations with the following properties, which we
justify as we state them.
\begin{enumerate}
\item There are $m$ input nodes labeled by the $m$ entries of $\bar{A}$,
$a_{1,j+1}$ through $a_{m,j+1}$. 
We call these $\bar{A}$-input nodes for short.
\item There are at least $mj-j(j+1)/2$ input nodes labeled by parameters
representing $\bar{Q}$, since this many parameters are needed to
represent a member of the Stiefel manifold.
We call these $\bar{Q}$-input nodes for short.
\item There are two types of computation nodes, addition and multiplication.
In other words, we assume that we do not do divisions, square roots, etc. 
Since we are only doing matrix multiplication, this is reasonable.
We note that any divisions or square roots in the overall algorithm
may be done in order to compute the parameters represented $\bar{Q}$.
Omitting these from consideration only lowers our lower bound
(though not by much).
\item There are no branches in the algorithm. In other words, the
way an entry of $\bar{R}$ is computed does not depend on the numerical
values. This assumption reflects current algorithms, but could in fact 
be eliminated as explained later.
\item Since the computation nodes only do multiplication and addition, 
we may view the output of each node as a polynomial in entries of $\bar{A}$
and parameters representing $\bar{Q}$.
\item We further restrict the operations performed so that the output of
any node must be a homogeneous linear polynomial in the entries of $\bar{A}$.
In other words, we never multiply two quantities depending on entries
of $\bar{A}$ to get a quadratic or higher order polynomial,
or add a constant or parameter depending on $\bar{Q}$ to an entry of
$\bar{A}$.  This is
natural, since the ultimate output is linear and homogeneous in $\bar{A}$, 
and any higher degree polynomial terms or constant terms would have to 
be canceled away.  No current or foreseeable algorithm (even Strassen based) 
would do this, and numerical stability would likely be lost.
\item There are $j$ output nodes labeled by the entries of $\bar{R}$,
$r_{1,j+1}$ through $r_{j,j+1}$.
\end{enumerate}

The final requirement means that multiplication nodes are only allowed
to multiply $\bar{Q}$-input nodes and homogeneous linear functions
of $\bar{A}$, including $\bar{A}$-input nodes.
Addition nodes may add homogeneous linear functions of $\bar{A}$ 
(again including $\bar{A}$-input nodes), but not add $\bar{Q}$-input nodes
to homogeneous linear functions of $\bar{A}$.
We exclude the possibility of adding or multiplying $\bar{Q}$-input nodes,
since the results of these could just be represented as additional
$\bar{Q}$-input nodes.

Thus we see that the algorithm represented by the DAG just described
outputs $j$ polynomials that are homogeneous and linear in $\bar{A}$.
Let $M$ be the total number of multiplication nodes in the DAG.
We now want to argue that unless $M \geq mj-j(j+1)/2$,
these output polynomials cannot possibly compute the right answer. 
We will do this by arguing that the dimension of
a certain algebraic variety they define is both bounded above by $M$,
and the dimension must be at least $mj-j(j+1)/2$ to get the right answer.

Number the output nodes from $1$ to $j$.
The output polynomial representing node $i$ can be written as
$\sum_{k=1}^m p_{k,i} (\bar{Q}) a_{k,j+1}$, where $p_{k,i}(\bar{Q})$ is
a polynomial in the values of the $\bar{Q}$-input nodes. According
to our rules for DAGs above, only multiplication nodes can introduce
a dependence on a previously unused $\bar{Q}$-input node, so
all the $p_{k,i}(\bar{Q})$ 
can only depend on $M$ independent parameters.

% I think the following argument is correct, but I should run it by
% an algebraic geometer. 
Finally, viewing each output node as a vector
of $m$ coefficient polynomials 
\linebreak
$(p_{1,i} (\bar{Q}),...,p_{m,i} (\bar{Q}))$,
we can view the entire output as a vector of $mj$ coefficient polynomials
$V(\bar{Q}) = (p_{1,1}(\bar{Q}),...,p_{m,j}(\bar{Q}))$, 
depending on $M$ independent parameters. 
This vector of length $mj$ needs to represent the set of
all $m$-by-$j$ orthogonal matrices. But the Stiefel manifold of
such orthogonal matrices has dimension $mj-j(j+1)/2$, so the surface
defined by $V$ has to have at least this dimension, i.e. $M \geq mj-j(j+1)/2$.

As an extension, we could add branches to our algorithm by noting that the
output of our algorithm would be piecewise polynomials, on regions
whose boundaries are themselves defined by varieties in the
same homogeneous linear polynomials. We can apply the above argument
on all the regions with nonempty interiors to argue that the same
number of multiplications is needed.

In summary, we have proven

\lemma{
Suppose we are doing the QR factorization of an $m$-by-$n$
matrix using any ``columnwise'' algorithm in the sense described
above. Then at least $mj - j(j+1)/2$ multiplications are required
to compute column $j+1$ of $R$, and at least
$\frac{mn^2}{4} - \frac{n^2}{8}(\frac{n}{2} + 1)$
multiplications to compute columns $\frac{n}{2}+1$ through $n$
of $R$.
}
\label{lemma:F_lowerbound}
\rm

\section{Related work}\label{S:related-work}

The central idea in this paper is factoring tall skinny matrices using
a tree-based Householder QR algorithm.  A number of authors previously
figured out the special case of a binary reduction tree for parallel
QR.  As far as we know, Golub et al.\ \cite{golub1988parallel} were the
first to suggest it, but their formulation requires $n \log P$
messages for QR of an $m \times n$ matrix on $P$ processors.  Pothen
and Raghavan \cite{pothen1989distributed} were the first, as far as we
can tell, to implement parallel TSQR using only $\log P$ messages.  Da
Cunha et al.\ \cite{cunha2002new} independently rediscovered parallel
TSQR.  

Other authors have worked out variations of the algorithm we call
``sequential TSQR''
\cite{buttari2007class,buttari2007parallel,gunter2005parallel,kurzak2008qr,quintana-orti2008scheduling,rabani2001outcore}.
They do not use it by itself, but rather as the panel factorization
step in the QR decomposition of general matrices.  The references
\cite{buttari2007class,buttari2007parallel,gunter2005parallel,kurzak2008qr,quintana-orti2008scheduling}
refer to the latter algorithm as ``tiled QR,'' which is the same as
our sequential CAQR with square blocks.  However, they use it in
parallel on shared-memory platforms, especially single-socket
multicore.  They do this by exploiting the parallelism implicit in the
directed acyclic graph of tasks.  Often they use dynamic task
scheduling, which we could use but do not discuss in this paper.
Since the cost of communication in the single-socket multicore regime
is low, these authors are less concerned than we are about minimizing
latency; thus, they are not concerned about the latency bottleneck in
the panel factorization, which motivates our parallel CAQR algorithm.
We also model and analyze communication costs in more detail than
previous authors did.

Here are recent examples of related work on sequential CAQR.
Gunter and van de Geijn develop a parallel out-of-DRAM QR
factorization algorithm that uses a flat tree for the panel
factorizations \cite{gunter2005parallel}.  Buttari et al.\ suggest
using a QR factorization of this type to improve performance of
parallel QR on commodity multicore processors \cite{buttari2007class}.
Quintana-Orti et al.\ develop two variations on block QR factorization
algorithms, and use them with a dynamic task scheduling system to
parallelize the QR factorization on shared-memory machines
\cite{quintana-orti2008scheduling}.  Kurzak and Dongarra use similar
algorithms, but with static task scheduling, to parallelize the QR
factorization on Cell processors \cite{kurzak2008qr}.

As far as we know, parallel CAQR is novel.  Nevertheless, there is a
body of work on theoretical bounds on exploitable parallelism in QR
factorizations.  These bounds apply to both parallel TSQR and parallel
CAQR if one replaces ``matrix element'' in the authors' work with
``block'' in ours.  Cosnard, Muller, and Robert proved lower bounds on
the critical path length $Opt(m,n)$ of any parallel QR algorithm of an
$m \times n$ matrix based on Givens rotations \cite{cosnard86}; it is
believed that these apply to any QR factorization based on Householder
or Givens rotations.  Leoncini et al.\ show that any QR factorization
based on Householder reductions or Givens rotations is P-complete
\cite{leoncini1999parallel}.  The only known QR factorization
algorithm in arithmetic NC (see \cite{csanky1976fast}) is numerically
highly unstable \cite{demmel1992trading}, and no work suggests that a
stable arithmetic NC algorithm exists.

Hong and Kung \cite{hong1981io} and Irony, Toledo, and Tiskin
\cite{irony2004communication} proved lower bounds on communication for
sequential and parallel matrix multiplication.  We are, as far as we
know, the first to attempt extending these bounds to LU and QR factorization.
Elmroth and Gustavson proposed a recursive QR factorization (see
\cite{elmroth1998new,elmroth2000applying}) which can also take
advantage of memory hierarchies.  It is future work to analyze whether
their algorithm satisfies the same lower bounds on communication as
does sequential CAQR.  It is natural to ask to how much of dense
linear algebra one can extend the results of this paper, that is
finding algorithms that attain communication lower bounds.  
For parallel LU with pivoting, see the technical report by
Grigori, Demmel, and Xiang \cite{grigori2008calu}, and for
sequential LU, see \cite{toledo1997locality}.

Block iterative methods frequently compute the QR factorization of a
tall and skinny dense matrix.  This includes algorithms for solving
linear systems $Ax = B$ with multiple right-hand sides (such as
variants of GMRES, QMR, or CG
\cite{vital:phdthesis:90,Freund:1997:BQA,oleary:80}), as well as block
iterative eigensolvers (for a summary of such methods, see
\cite{templatesEigenBai,templatesEigenLehoucq}).  In practice,
modified Gram-Schmidt orthogonalization is usually used when a 
(reasonably) stable
QR factorization is desired.  Sometimes unstable methods (such as
CholeskyQR) are used when performance considerations outweigh
stability.  Eigenvalue computation is particularly sensitive to the
accuracy of the orthogonalization; two recent papers suggest that
large-scale eigenvalue applications require a stable QR factorization
\cite{lehoucqORTH,andrewORTH}.  Many block iterative methods have
widely used implementations, on which a large community of scientists
and engineers depends for their computational tasks.  Examples include
TRLAN (Thick Restart Lanczos), BLZPACK (Block Lanczos), Anasazi
(various block methods), and PRIMME (block Jacobi-Davidson methods)
\cite{TRLANwebpage,BLZPACKwebpage,BLOPEXwebpage,irbleigs,TRILINOSwebpage,PRIMMEwebpage}.

\section{Conclusions and Open Problems}
\label{sec:Conclusions_optimal}

We have shown that known bandwidth lower bounds
for parallel and sequential $\Theta(n^3)$ matrix multiplication
imply latency lower bounds,
shown such bounds apply to both LU and QR algorithms,
presented some new and some old QR algorithms
that attain these bounds, and referred to LU
algorithms in the literature that attain at least
some of these bounds. Whether a sequential LU
algorithm exists attaining the latency lower
bound is an open question.

There are numerous ways in which one could hope to
extend these results. One natural conjecture is that
the bounds apply to other $\Theta(n^3)$ dense linear algebra
routines, such as eigenvalue problems, and if they do,
we would want to find algorithms that attain them. Another question
is finding analogous communication lower bounds for 
asymptotically faster dense linear algebra algorithms 
like thosed based on Strassen's algorithm, 
or indeed of any matrix multiplication algorithm, 
based on Raz's theorem converting any matrix
multiplication algorithm to be ``Strassen-like''
(bilinear noncommutative) \cite{raz2003complexity}.

But the following question is of more practical importance.
Our TSQR and CAQR algorithms have been described and analyzed 
in most detail for simple machine models: either
sequential with two levels of memory hierarchy (fast and slow),
or a homogeneous parallel machine, where each processor is
itself sequential. Real computers are more complicated, with
many levels of memory hierarchy and many levels of parallelism
(multicore, multisocket, multinode, multirack, \dots) all with
different bandwidths and latencies. So it is natural to ask
whether our algorithms and optimality proofs can be extended
to these more general situations. We hinted at
how TSQR could be extended to general
reduction trees in Section~\ref{sec:TSQR_optimal}, which
could in turn be chosen depending on the architecture.
But we have not discussed CAQR, which we do here.

We again look at the simpler case of matrix multiplication
for inspiration. Consider the sequential case, with
$k$ levels of memory hierarchy instead of 2, where
level 1 is fastest and smallest with $W_1$ words of memory,
level 2 is slower and larger with $W_2$ words of memory,
and so on, with level $k$ being slowest and large enough
to hold all the data. By dividing this hierarchy into
two pieces, levels $k$ through $i+1$ ("slow")
and $i$ through 1 ("fast"), we can apply the theory in
Section~\ref{SS:MMlowerbounds} to get lower bounds
on bandwidth and latency for moving data between levels $i$
and $i+1$ of memory. So our goal expands to finding
a matrix multiplication algorithm that attains not just
1 set of lower bounds, but $k-1$ sets of lower bounds,
one for each level of the hierarchy. 

Fortunately, as is well known, the standard approach to tiling 
matrix multiplication achieves all these lower bounds simultaneously, 
by simply applying it recursively: level $i+1$ holds submatrices
of dimension $\Theta(\sqrt{W_{i+1}})$, and multiplies them by tiling
them into submatrices of dimension $\Theta(\sqrt{W_i})$, and so on.

The analogous observation is true of parallel matrix multiplication
on a hierarchical parallel processor where each node in the parallel 
processor is itself a parallel processor (multicore, multisocket, 
multirack, \dots).

We believe that this same recursive hierarchical approach applies
to CAQR (and indeed much of linear algebra) but there is a catch:
Simple recursion does not work, because the subtasks are not all
simply smaller QR decompositions. Rather they are a mixture of
tasks, including smaller QR decompositions and operations like 
matrix multiplication. Therefore we still expect that the same
hierarchical approach will work: if a subtask is matrix multiplication
then it will be broken into smaller matrix multiplications as
described above, and if it is QR decomposition, it will be broken into
smaller QR decompositions and matrix multiplications.

There are various obstacles to this simple approach.
First, the small QR decompositions generally have structure,
e.g., a pair of triangles. To exploit this structure fully
would complicate the recursive decomposition. (Or we could 
ignore this structure, perhaps only on the smaller
subproblems, where the overhead would dominate.)

Second, it suggests that the data structure with which the
matrix is stored should be hierarchical as well, with
matrices stored as subblocks of subblocks \cite{elmroth2004recursive}.
This is certainly possible, but it differs significantly
from the usual data structures to which users are accustomed.
It also suggests that recent approaches based on decomposing
dense linear algebra operations into DAGs of subtasks \cite{buttari2007class,boboulin2008issues,kurzak2008qr,quintana-orti2008scheduling,quintana-orti2008design}
may need to be hierarchical, rather than have a single layer
of tasks. A single layer is a good match for the single socket 
multicore architectures that motivate these systems, but may
not scale well to, e.g., petascale architectures.

Third, it is not clear whether this approach best
accommodates machines that mix hierarchies of parallelism
and memory. For example, a multicore / multisocket / multirack 
computer will have also have disk, DRAM and various caches,
and it remains to be seen whether straightforward recursion
will minimize bandwidth and latency everywhere that
communication takes place within such an architecture.

Fourth and finally, all our analysis has assumed homogeneous machines,
with the same flop rate, bandwidth and latency in all components. This
assumption can be violated in many ways, for example, by asymmetric
read and write bandwidths, by having different bandwidth and latency
between racks, sockets, and cores on a single chip, or by having some
specialized floating point units like GPUs.

It is most likely that an adaptive, ``autotuning'' approach
will be needed to deal with some of these issues, just
as it has been used for the simpler case of a matrix
multiplication.  Addressing all these issues is future work.

\newpage
\bibliographystyle{siam}
\bibliography{qr}
\end{document}